\documentstyle[12pt]{article}
  \def\Bbb{\bf}
  
  \def\boxtimes{\Box\hspace{-0.76em}\times}
 \catcode`@=11
 \def\@date{September, 1994}
 \catcode`@=12
\newcommand{\sgn}{\mathop{\rm sgn}}
\newcommand{\Row}{\mathop{\rm Row}}
\newcommand{\Col}{\mathop{\rm Col}}

\newcommand{\IM}{\mathop{\rm Im}}
\newcommand{\Sym}{\mathop{\rm Sym}\mbox{}}
\newcommand{\Poly}{\mathop{\rm Poly}\mbox{}}
\newcommand{\St}{\mathop{\rm St}}
\newcommand{\Gr}{\mathop{\rm Gr}}
\newcommand{\id}{\mathop{\rm id}}

\newcommand{\rest}{\mathop{\rm rest}\mbox{}}
\newcommand{\rank}{\mathop{\rm rank}}
\newcommand{\pr}{ {\mathop{\rm pr}} }
\newcommand{\wt}{ {\mathop{\rm wt}} }
\newcommand{\tr}{\mathop{\rm tr}}
\newcommand{\gap}{ {\mathop{\rm gap}} }
\newcommand{\Hom}{\mathop{\rm Hom}}
\newcommand{\Char}{ {\mathop{\rm char}} }
\newcommand{\Span}{ {\mathop{\rm Span}} }
\newcommand{\diag}{\mathop{\rm diag}}
\newtheorem{thm}{Theorem}
\newtheorem{lem}[thm]{Lemma}

\newtheorem{prop}[thm]{Proposition}
\newtheorem{conj}[thm]{Conjecture}
\newcommand{\eqdef}{\stackrel{\rm def}{=}}
\newcommand{\C}{{\Bbb C}}
\newcommand{\R}{{\Bbb R}}
\newcommand{\Z}{{\Bbb Z}}
\newcommand{\NN}{{\Bbb N}}
\newcommand{\PP}{{\Bbb P}}
\newcommand{\al}{{\alpha}}

\newcommand{\alD}{{\alpha_D}}
\newcommand{\beD}{{\beta_D}}

\newcommand{\lam}{ \lambda }
\newcommand{\FF}{ { \cal F } }
\newcommand{\FD}{ \FF_D }
\newcommand{\Fgen}{ \FF^{\mbox{\rm \tiny gen}} }
\newcommand{\II}{ { \cal I } }
\newcommand{\ID}{ \II_D }
\newcommand{\LL}{ {\cal L} }
\newcommand{\LD}{ \LL_D }
\newcommand{\EE}{ {\cal E} }
\newcommand{\OO}{ {\cal O} }
\newcommand{\DV}{ \Delta^D V}

\newcommand{\llx}{ \stackrel{\rm lex}{<}  }
\newcommand{\glx}{ \stackrel{\rm lex}{>}  }
\newcommand{\leqlx}{ \stackrel{\rm lex}{\leq}  }
\newcommand{\linit}{ \stackrel{\rm init}{\subset} }
\newcommand{\ginit}{ \stackrel{\rm init}{\supset} }
\newcommand{\Dh}{ \widehat{D} }
\newcommand{\Ch}{ \widehat{C} }

\newcommand{\tB}{ \stackrel{B}{\times} }
\newcommand{\subneq}{ \stackrel{\neq}{\subset} }

\title{A Borel-Weil Construction for Schur Modules}
\author{Peter Magyar}

\begin{document}

\maketitle

\begin{center} {\bf Abstract} \end{center}

{\small
\noindent
We present a generalization of the classical
Schur modules of $GL(N)$ exhibiting the
same interplay among algebra, geometry, and combinatorics.
A generalized Young diagram $D$
is an arbitrary finite subset of $\NN \times \NN$.
For each $D$, we define
the Schur module $S_D$ of $GL(N)$.
We introduce a projective variety $\FF_D$
and a line bundle $\LL_D$,
and describe the Schur module
in terms of sections of $\LL_D$.

For diagrams with the ``northeast'' property,
 $$
 (i_1,j_1),\ (i_2, j_2) \in D  \Rightarrow
(\min(i_1,i_2),\max(j_1,j_2)) \in D ,
$$
which includes the skew diagrams,
we resolve the singularities of $\FD$
and show analogs of Bott's and Kempf's vanishing theorems.
Finally, we apply the Atiyah-Bott Fixed Point Theorem
to establish a Weyl-type character formula
of the form:
$$
{\Char}_{S_D}(x) =
\sum_t {x^{\wt(t)} \over \prod_{i,j} (1-x_i x_j^{-1})^{d_{ij}(t)}} \ ,
$$
where $t$ runs over certain standard tableaux of $D$.

Our results are valid over fields of arbitrary characteristic.
}
\\[1em]
{\large \bf Introduction}
\\[1em]
The two main branches of the representation theory of the general
linear groups $G = GL(N,F)$ began with
the geometric Borel-Weil-Bott theory
and the combinatorial analysis of Schur, Young, and Weyl.
In the case when $F$ is of characteristic zero,
the geometric theory  realizes the irreducible representation
of $G$ with highest weight $\lam$ as the sections of a
line bundle $\LL_{\lam}$ over the flag variety $\FF = G/B$.
By contrast, the Schur-Weyl construction obtains
this representation as a subspace of $V^{\otimes k}$,
where $V = F^N$,
by a process of symmetrization and anti-symmetrization defined
by $\lam$ considered as a Young diagram.
(See~\cite{FH} for an accessible reference.)

Combinatorists have examined other symmetrization operations
on $V^{\otimes k}$, such as those associated to skew
diagrams, and recently more general diagrams $D$
of squares in the plane
(~\cite{ABW},~\cite{JP},~\cite{Kr},~\cite{KP},~\cite{MP},
{}~\cite{RS1},~\cite{RS2},~\cite{RS3},
{}~\cite{W}).
We call the resulting $G$-representations
the Schur modules $S_D \subset V^{\otimes k}$.
(In characteristic zero, $S_D$ is irreducible exactly
when $D$ is a Young diagram.)
Kraskiewicz and Pragacz~\cite{KP} have shown
that the characters of $S_D$, for $D$ running through the
inversion diagrams of the symmetric group on $N$ letters,
give an algebraic description of the Schubert calculus
for the cohomology of the flag variety $\FF$.
(More precisely, the Schubert polynomials are characters
of {\em flagged} Schur modules. We deal with this
case in~\cite{M1}.)

In this paper, we attempt to combine the combinatorial
and the geometric approaches.
We give a geometric definition (valid for all characteristics)
for the $G$-module $S_D$.  That is, for any finite
set $D \subset \NN \times \NN$, we produce $S_D$
as the space of sections of a line bundle over a
projective variety $\FF_D$,
the configuration variety of $D$.
(This is proved only
for diagrams with a ``direction'' property,
but a weaker statement is shown
for general diagrams.)
Our picture reduces to that of
Borel-Weil when D is a Young diagram.
See also~\cite{BD1},~\cite{BD2},
where similar varieties are introduced.
We prove a conjecture of
V. Reiner and M. Shimozono
asserting the duality between
the Schur modules of two diagrams whose
disjoint union is a rectangular diagram.

We can carry out a more detailed analysis
for diagrams satisfying a
direction condition such as the
northeast condition
$$
 (i_1,j_1),\ (i_2, j_2) \in D  \Rightarrow
(\min(i_1,i_2),\max(j_1,j_2)) \in D .
$$
To accord with the literature, we will deal
exclusively with north{\em west} diagrams,
but since the modules and varieties with which we
are concerned do not change
(up to isomorphism)
if we switch one row of the diagram with another or
one column with another, everything we will say
applies with trivial modifications
to skew, inversion, Rothe,
and column-convex diagrams, and
diagrams satisfying any direction condition (NE, NW, SE, SW).

In this case,
we find an explicit resolution
of singularities of $\FD$,
and we use Frobenius splitting arguments of Wilberd
van der Kallen (based on work of Mathieu,
Polo, Ramanathan, et al.)
to show the vanishing of certain higher
cohomology groups.
In particular, the configuration varieties
are projectively normal and have rational singularities.
This allows us to
apply the Atiyah-Bott Fixed Point Theorem
to compute the character of the Schur modules.

For more general diagrams, the above program breaks down
because we lack a suitable desingularization of $\FD$.
We can carry it through, however, for diagrams with
at most 3 rows, for which the configuration variety is
the space of triangles~\cite{FM},~\cite{M2}.

Those interested only in the algebraic and combinatorial side
of our results can find the definitions and statements in
Section \ref{Schur modules and Weyl modules},
Theorem \ref{complement thm}
of Section \ref{complementary diagrams},
and Section \ref{character formula}.
Our discussion of geometry begins with Section \ref{configuration
varieties}.

\vspace{1em}

\noindent
{\small
{\bf Contents.}\
{\bf 1 $GL(N)$ Modules}\
1.1 Schur modules \
1.2 Weyl modules \
{\bf 2 Configuration varieties} \
2.1 Definitions and examples \
2.2 Diagrams with at most $N$ rows \
2.3 Complementary diagrams \
{\bf 3 Resolution of singularities} \
3.1 Northwest and lexicographic diagrams \
3.2 Blowup diagrams \
3.3 Blowup varieties \
3.4 Intersection varieties \
3.5 Smoothness and equations defining varieties \
{\bf 4 Cohomology of line bundles}\
4.1 Frobenius splittings of flag varieties \
4.2 Frobenius splittings of Grassmannians \
4.3 Monotone sequences of permutations \
{\bf 5  A Weyl character formula}\
5.1 Fixed points and tangent spaces \
5.2 The character formula \
5.3 Betti numbers \ \
References
}

\vspace{2em}

\noindent
{\bf Acknowledgements.}
I am deeply indebted
to Wilberd van der Kallen,
who supplied the essential Propositions
\ref{splitting}
and \ref{rational sing},
and gave many valuable suggestions.
I also wish to thank
Victor Reiner and Mark Shimozono for making me acquainted with
their conjectures, as well as Eduard Looijenga, Tony Springer, and
S.P. Inamdar
for helpful discussions.

\section{$GL(N)$ modules}
\label{Schur modules and Weyl modules}

\subsection{Schur modules}

Given a finite set $T$,
we will also use the symbol $T$ to denote
the order $|T|$ when appropriate.
Thus $GL(T) \eqdef GL(|T|)$, etc.
Let $\Sigma_T$ be the symmetric
group permuting the elements of $T$.
For any left $G$-space $X$, $\Sigma_T$ acts on the
right, and $G$ acts on the left, of the cartesian
product $X^T$ by:
$$
g(x_{t_1}, x_{t_2}, \ldots ) \pi =
(g x_{\pi t_1}, g x_{\pi t_2}, \ldots )  .
$$
A {\em diagram} is a finite subset of
$\NN \times \NN$.
Its elements $(i,j) \in D$ are called {\em squares}.
We shall often think of $D$ as a sequence
$(C_1,C_2,\ldots,C_r)$ of columns $C_j \subset \NN$.
The Young diagram corresponding to
$\lam = (\lam_1 \geq \lam_2 \geq \cdots \geq \lam_N \geq 0)$ is the
set $\{(i,j) \mid 1\leq j\leq N,\  1\leq i\leq \lam_j \}$.
For any diagram $D$, we let
$$
\Col(D) = \{\pi \in \Sigma_D \mid \pi(i,j) = (i',j) \ \exists i'\}
$$
be the group permuting the squares of $D$
within each column, and we define $\Row(D)$ similarly
for rows.

Let $F$ be a field.
We shall always write $G = GL(N,F)$, $B = $ the subgroup
of upper triangular matrices, $H =$ the subgroup of diagonal
matrices, and $V = F^N$ the
defining representation.

Now let $F$ have characteristic zero.
Define the idempotents
$\al_D$, $\beta_D$ in the group algebra
$F[\Sigma_D]$ by
$$
\alD = {1 \over |\Row D|} \sum_{\pi \in \Row D} \pi, \ \ \ \
\beD = {1 \over |\Col D|} \sum_{\pi \in \Col D} \sgn(\pi) \pi ,
$$
where $\sgn(\pi)$ is the sign of the permutation.
Define the {\em Schur module}
$$
S_D \eqdef V^{\otimes D} \alD \beD \subset V^{\otimes D},
$$
a representation of $G$.

Note that we get
an isomorphic Schur module
if we change the diagram by permuting the rows or the
columns (i.e., for some permutation $\pi: \NN \rightarrow \NN$,
changing $D = \{(i,j)\}$ to
$D' = \{(\pi(i), j) \mid (i,j) \in D\}$, and similarly for columns).

\subsection{Weyl modules}

Let $W = V^*$, the dual of the defining representation of
$G=GL(N,F)$, where $F$ is an infinite field.
Given a diagram $D$,
define the alternating product
{\em with respect to the columns}
$$
{\bigwedge} ^D W = \{ f:V^D \rightarrow F \mid
f \mbox{ multilinear, and } f(v \pi) = \sgn(\pi) f(v) \ \forall \pi
\in \Col(D) \},
$$
where {\em multilinear} means $f(v_1,\ldots,v_d)$ is $F$-linear
in each of the $d = |D|$ variables.
Consider the multidiagonal
{\em with respect to the rows}
$$
\DV = \Delta^{R_1} V \times
\Delta^{R_2} V \times \cdots \subset V^{R_1} \times V^{R_2}
\times \cdots = V^D ,
$$
where $R_1, R_2, \ldots $ are the rows of $D$.
Now define the {\em Weyl module}
$$
W_D \eqdef {\bigwedge} ^D W \mid_{\DV},
$$
where $\mid_{\DV}$ denotes restriction of functions from
$V^D$ to $\DV$.
Since $\DV$ is stable under the diagonal action of
$G$, $W_D$ is naturally a $G$-module.

\vspace{1em}

\noindent
{\bf Remark.} For $F$ a finite field,
we make the following modification.
Consider $W = W(F) \hookrightarrow W(\bar{F})$, where
$\bar{F}$ is the algebraic closure.
That is, identify
$$
W = \{ f:\bar{F}^N \rightarrow \bar{F} \mid
f \mbox{ is $\bar{F}$-linear, and }
f(F^N) \subset F \}.
$$
Then define
$$
W_D \eqdef {\bigwedge} ^D W \mid_{\DV(\bar{F})}.
$$
This keeps the restriction map from killing
nonzero tensors which happen to vanish on the finite set
$\DV(F)$.

With this definition, $W_D$ clearly has the base change property
$W_D(L) = W_D(F) \otimes_F L$ for any extension of
fields $F \subset L$.

Now consider $W_D(\Z)$.  This is a free $\Z$-module, since it is
a submodule of the $\Z$-valued functions on $\DV$.
Suppose $D$ satisfies a direction condition.  Then our vanishing
results of Proposition \ref{config split} (a),
along with the appropriate universal
coefficient theorems, can be used to show that for any field $F$,
$$
W_D(F) = W_D(\Z) \otimes_{\Z} F .
$$

\begin{prop}
If $F $ has characteristic zero, then $W_D \cong
S_D^*$ as $G$-modules.
\end{prop}

\noindent
{\bf Proof.} $S_D$ is the image of the composite mapping
$$
V^{\otimes D} \alD \hookrightarrow
V^{\otimes D} \stackrel{\beD}{\rightarrow}
V^{\otimes D} \beD .
$$
For $W = V^*$, write
$$
W^{\otimes D} = \{f: V^D \rightarrow F \mid f \mbox{ multilinear}\},
$$
$$
\Sym ^D W =
\{f: V^D \rightarrow F \mid f \mbox{ multilinear, and }
f(v\pi) = f(v) \ \forall \pi \in \Row(D)\}.
$$
Now, representations of $F[\Sigma_D]$ are completely reducible,
so $S_D^*$ is the image of
$$
W^{\otimes D} \beD \hookrightarrow
W^{\otimes D} \stackrel{\alD}{\rightarrow}
W^{\otimes D} \alD ,
$$
and $ W^{\otimes D} \beD \cong \bigwedge^D W$,\
$W^{\otimes D} \alD \cong \Sym^D W$.

Now, let
$$
{\Poly}^D W =
\{f: V^l \rightarrow F \mid f \mbox{ homog poly of multidegree }
(R_1, \ldots, R_l) \},
$$
where $l$ is the number of rows of $D$.
Then we have a $G$-equivariant map
$$
\rest_{\Delta} : \Sym  ^D W \rightarrow \Poly ^D W
$$
restricting functions from $V^D$ to the row-multidiagonal $\DV \cong V^l$.
It is well known that $\rest_{\Delta}$ is an isomorphism:
the symmetric part of a tensor algebra is isomorphic
to a polynomial algebra.

Thus we have the commutative diagram
\begin{eqnarray*}
{\bigwedge} ^D  W & \hookrightarrow
W^{\otimes D} \stackrel{\alD}{\rightarrow}
& \Sym  ^D W \\
|| & || & \downarrow \rest_{\Delta} \\
{\bigwedge} ^D W & \hookrightarrow
W^{\otimes D} \stackrel{\alD}{\rightarrow}
& \Poly ^D W  .
\end{eqnarray*}
Now, the image in the top row is $S_D^*$,
the image in the bottom row
is $W_D$, and all the vertical maps are isomorphisms,
so we have $\rest_{\Delta} : S_D^* \tilde{\rightarrow} W_D$
an isomorphism. $\bullet$


If $D = \lam$ a Young diagram,
then $W_D$ is  isomorphic to Carter and Lusztig's
dual Weyl module for $G = GL(N,F)$.
This will follow from Proposition \ref{Young diagram} in the
following section.

\section{Configuration varieties}
\label{configuration varieties}

{\bf N.B.} Although our constructions
remain valid over $\Z$, for simplicity we will
assume for the remainder of this paper
that
{\em $F$ is an algebraically closed field}.

\subsection{Definitions and examples}
\label{definitions and examples}

Given a finite set C (a column), and $V=F^N$, consider
$ V^C \cong M_{N\times C}(F)$, the $N \times |C|$ matrices,
with a right multiplication of  $GL(C)$.
Let
$$
\St(C) = \{X \in V^C \mid \rank X = |C| \},
$$
the Stiefel manifold, and
$$
\Gr(C) = \St(C)/GL(C),
$$
the Grassmannian.

Also, let
$$
{\cal L}_C = \St(C) \stackrel{GL(C)}{\times} {\det} ^{-1} \rightarrow
\Gr(C)
$$
be the Plucker determinant bundle,
whose sections are regular functions $f: \St(C) \rightarrow F$ with
$f(XA) = \det(A) f(X) \ \forall A \in GL(C)$.
In fact, such global sections can be extended to polynomial
functions
$f: V^C \rightarrow F$.

For a diagram $D$ with columns $C_1, C_2,\ldots$,
we let
$$
\St(D) = \St(C_1) \times \St(C_2) \times \cdots,
\mbox{  }
\Gr(D) = \Gr(C_1) \times \Gr(C_2) \times \cdots,
\mbox{  }
\LD = {\cal L}_{C_1} \boxtimes {\cal L}_{C_2} \boxtimes \cdots.
$$

Recall that $\DV \subset V^D$ is the {\em row} multidiagonal
(as opposed to the column constructions above).
Let
$$
\FD^o \eqdef \IM\left[\DV \cap \St(D) \rightarrow \Gr(D)\right],
$$
and define the {\em configuration variety of $D$} by
$$
\FD = \overline{\FD^o} \subset \Gr(D),
$$
the Zariski closure of $\FD^o$ in $\Gr(D)$.
We denote the restriction of $\LD$ from $\Gr(D)$
to $\FD$ by the same symbol $\LD$.

Some properties follow immediately from the definitions.
For instance, $\FD$ is an irreducible variety.
Just as for Schur modules
and Weyl modules,
changing the diagram by permuting the rows
or the columns gives an isomorphic configuration variety and
line bundle.
If we add a column $C$ to $D$ which already
appears in $D$, we get an isomorphic configuration
variety, but the line bundle
is twisted to have higher degree.
Since $\LD$ gives the Plucker embedding on $\Gr(D)$,
it is very ample on $\FD$.

\vspace{1em}

\noindent
{\bf Examples. } Set $N=4$.  Consider the diagrams
$$
D_1 = \begin{array}{ccc}
 \Box & & \\
 \Box & \Box & \\
      & \Box &
\end{array}
\ \ \
D_2 = \begin{array}{ccc}
 \Box &      & \\
 \Box & \Box & \Box \\
      &      & \Box
\end{array}
\ \ \
D_3 = \begin{array}{ccc}
 \Box & \Box & \\
      & \Box &      \\
      & \Box & \Box
\end{array}
$$

Identifying $\Gr(k,F^N)$ with $\Gr(k-1,\PP^{N-1}_F)$, we may
consider the $\FD$'s as varieties of configurations in $\PP^3$:\\
(1)  $\FF_{D_1}$ is the variety of pairs $(l,l')$, where
$l,l'$ are intersecting lines in $\PP^3$.
It is singular at the locus where the two lines coincide.\\
(2)  $\FF_{D_2}$ is the variety of triples $(l, p, l')$
of two lines and a point which lies on both of them.
The variety is smooth:  indeed, it is a fiber bundle over the partial
flag variety of a line containing a point.  There is an obvious
map $\FF_{D_2} \rightarrow \FF_{D_1}$, which is birational,
and is in fact a small resolution of singularities.
(C.f. Proposition \ref{smooth}.)\\
(3)  $\FF_{D_3}$ is the variety of planes with two marked
points (which may coincide).  $\FF_{D_3}^o$ is the locus where the
marked points are distinct. The variety is smooth as in
the previous example.  \\
$$
D_4 = \begin{array}{ccc}
 \Box & \Box & \\
 \Box &      & \Box \\
      & \Box & \Box
\end{array}
\ \ \
D_5 = \begin{array}{ccc}
 \Box & \Box & \Box \\
 \Box &      &      \\
      & \Box &      \\
      &      & \Box
\end{array}
D_6 = \begin{array}{cccc}
\Box  &      &      &     \\
      & \Box &      &     \\
      &      & \Box &     \\
      &      &      &   \Box
\end{array}
\ \ \
$$
(4) $\FF_{D_4}$ is the variety of triples of coplanar lines.\\
(5) $\FF_{D_5}$ is the variety of triples of lines with a common point.
This is the
projective dual of the previous variety,
since the diagrams are complementary
within a $4 \times 3$ rectangle (up to permutation of rows and
columns).
(See Theorem \ref{complement thm}.)
The variety of triples of lines which intersect
pairwise cannot be described by a single diagram, but consists of
$\FF_{D_4} \cup \FF_{D_5}$. (See Section \ref{intersection varieties}.)\\
(6) $\FF_{D_6} \cong (\PP^{3})^4$ contains the
$GL(N)$-invariant subvariety where all four points
in $\PP^3$ are colinear.  Since the cross-ratio
is an invariant of four points on a line,
this subvariety contains infinitely many $GL(N)$ orbits.
$$
D_7 = \begin{array}{cccc}
 \Box & \Box & \Box &      \\
      & \Box & \Box & \Box \\
      &      & \Box &      \\
      &      &      & \Box
\end{array}
\ \ \
D_8 = \begin{array}{ccccc}
 \Box & \Box & \Box &      &      \\
      & \Box & \Box & \Box & \Box \\
      &      & \Box &      &      \\
      &      &      & \Box &
\end{array}
$$
(7)  $\FF_{D_7} \cong G \tB X_{\lam}$, the $G$-orbit version of
the Schubert variety $X_{\lam} \subset \Gr(2,4)$
associated to the partition $\lam = (1,2)$.
This is the smallest example of a singular Schubert variety. \\
(8) $\FF_{D_8}$ is a smooth variety
which maps birationally
to $\FF_{D_7}$ by forgetting the point associated to
the last column.  In fact, this is essentially the same
resolution as (1) and (2) above.
Such resolutions of singularities can be given
for arbitrary Schubert varieties of $G = GL(N)$,
and generalize Zelevinsky's resolutions
in~\cite{Z}.  C.f. Section \ref{resolution}.
$\bullet$

\begin{thm}
\label{weyl mod is image}
If $F$ is an algebraically closed field, then
$$
W_D \cong
\IM\left[\rest_{\Delta}: H^0(\Gr(D), \LD) \rightarrow H^0(\FD,\LD)\right],
$$
where $\rest_{\Delta}$ is the restriction map.
\end{thm}

\noindent
{\bf Proof.}  Note that for $GL(D) = GL(C_1) \times GL(C_2) \times \cdots$,
$$
H^0(\Gr(D), \LD) =
\{ f:V^D \rightarrow F \mid
  f(XA) = \det(A) \, f(X) \ \forall A \in GL(D)\},
$$
and recall
$$
{\bigwedge} ^D W = \{ f:V^D \rightarrow F \mid
f \mbox{ multilinear, and } f(v \pi) = \sgn(\pi) \ f(v) \ \forall \pi
\in \Col(D) \}.
$$
But in fact these sets are equal, because a multilinear, anti-symmetric
function $g:V^C\rightarrow F$ always satisfies
$g(XA) = \det(A) \ g(X) \, \forall A \in GL(C)$.
Now $W_D$ and $H^0(\FD,\LD)$ are gotten by restricting
functions in these identical sets to $\DV$, so we are done. $\bullet$

\subsection{Diagrams with at most $N$ rows}

We say $D$ has $\leq N$ rows if $(i,j) \in D \Rightarrow 1\leq i\leq N$.
\begin{prop}
\label{bigorbit}
If $D$ has $\leq N$ rows, then $\FD$ has an open dense $GL(N)$-orbit
$\Fgen_D$.
\end{prop}

\noindent
{\bf Proof.}
Let $D$ have columns $C_1, C_2,\ldots$.
Consider a sequence of vectors $X =(v_1,\ldots,v_n) \in V^N$.
For $C = \{i_1, i_2, \ldots \} \subset \{1,\ldots,N\}$,
define $X(C) \eqdef \Span_F(v_{i_1},v_{i_2},\ldots ) \in \Gr(C)$
(for $X$ sufficiently general).
Consider an element $g \in GL(N)$ as a sequence of
column vectors $g = (v_1,\ldots,v_n)$.
Then
$$
g(C) = g\cdot \Span_F(e_{i_1},e_{i_2},\ldots )
= g\cdot (I(C_1), I(C_2), \ldots) ,
$$
where $e_i$ denotes the $i$-th coordinate vector and
$I$ the identity matrix.

Now define the map
$$
\begin{array}{rccc}
\Psi: & V^N               & \rightarrow   & \DV \subset V^D \\
      & (v_1,\ldots,v_N)  & \mapsto       & (v_i)_{(i,j)\in D},
\end{array}
$$
where $(u_{ij})_{(i,j)\in D} $ denotes an element of $V^D$.
Then the composite
$$
V^N \stackrel{\Psi}{\rightarrow}  \DV \rightarrow  \FD^o
$$
is an onto map taking
$g \mapsto (g(C_1), g(C_2),\ldots) = g\cdot (I(C_1), I(C_2), \ldots)$.
Since $GL(N)$ is dense in $V^N$,
its image is dense in $\FD^o$, and hence the composite
image $\Fgen_D \eqdef G \cdot (I(C_1),I(C_2),\ldots)$
is a dense $G$-orbit in $\FD$.
 $\bullet$

\begin{prop}
\label{Young diagram}
If $D$ is the Young diagram associated to a
dominant weight $\lam$ of $GL(N)$, then: \\
(a) $\FD \cong G/P$, a quotient of the flag variety $\FF = G/B$. \\
(b) The Borel-Weil line bundle
$\LL_{\lam} \eqdef G \stackrel{B}{\times} (\lam^{-1}) \rightarrow \FF$ is
the pullback of $\LD$ under the projection $\FF\rightarrow \FD$. \\
(c) $   \rest_{\Delta} :
	H^0(\Gr(D),\LD) \rightarrow H^0(\FD,\LD)$
is surjective,
and $W_D \cong H^0(\FD,\LD)$. \\
\end{prop}

\noindent
{\bf Proof.}  (a) Let $\mu = (\mu_1 \geq \mu_2 \geq \cdots) = \lam^t $,
the transposed diagram, and let $P = \{(x_{ij}) \in GL(N) \mid
x_{ij} = 0 \mbox{ if } \exists k , \ i > \mu_k \geq j > \mu_{k+1} \}$,
a parabolic subgroup of $G$.
Then $G/P$ is the space of partial flags
$V = F^N \supset V_1 \supset V_2 \supset \cdots $
consisting of subspaces $V_j$ with $dim(V_j) = \mu_j$.
Clearly $G/P \cong \FD$. \\
(b)  Let $\Psi: G \rightarrow \DV \cap \St(D) $ be the map in the proof
of the previous proposition. Then the map
\begin{eqnarray*}
G \stackrel{P}{\times} F_{\lam^{-1}} &  \rightarrow  &
\LD = (\St(D) \stackrel{GL(D)}{\times} {\det} _D^{-1}) \mid_{\FD} \\
(g, \alpha) & \mapsto &  (\Psi(g), \alpha)
\end{eqnarray*}
is a $G$-equivariant bundle isomorphism.  Then (b) follows
by standard arguments. \\
(c) The surjectivity is a special case of
Proposition \ref{config split}
in Section \ref{cohomology}.
(See also~\cite{In}.)
The other statement then follows by Prop \ref{weyl mod is image}.
$\bullet$

\subsection{Complementary diagrams}
\label{complementary diagrams}

\begin{thm}
\label{complement thm}
Suppose the rectangular diagram
Rect $= \{1,\ldots,N\} \times \{1,\ldots,r\}$ is
the disjoint  union of two diagrams $D$, $D'$.
Let $W_D$, $W_{D'}$ be the corresponding Weyl modules for $G =GL(N,F)$.
Then: \\
(a) there is an $F$-linear bijection
$\tau:W_D \rightarrow W_{D'} $
such that
$\tau(g w) = \det_{N\times N}^r(g')\  g'\ \tau(w)$,
where $g'$ is the inverse transpose in $GL(N)$ of the matrix $g$; \\
(b) the characters obey the relation
$\Char \, {W_{D'}}(h) = \det_{N\times N}^r(h)\   \Char \, {W_D} (h^{-1}) ,$
for diagonal matrices $h\in G$ ; \\
(c) if $F$ has characteristic zero, then as $G$-modules
$$
W_{D'} \cong {\det}\mbox{}^{-r} \otimes W_{D}^*  \ \ \ \mbox{ and }\ \ \
S_{D'} \cong {\det}\mbox{}^r \otimes S_{D}^* .
$$
\end{thm}

\noindent
{\bf Proof.}  (a) Given $C \subset \{1,\ldots,N\} $
(a column set), we considered above the Plucker line bundle
$$
\St(C) \stackrel{GL(C)}{\times} \det\mbox{} ^{-1} \rightarrow \Gr(D) .
$$
We may equally well write this as
$$
GL(N) \stackrel{P_C}{\times} \det\mbox{} _C^{-1} \rightarrow \Gr(D) ,
$$
where $P_C \eqdef \{ (x_{ij}) \in GL(N) \mid
x_{ij} = 0 \mbox{ if } i \not\in C, \ j \in C\}$
is a maximal parabolic subgroup of $GL(N)$ (not necessarily
containing $B$), and
$\det_C : P_C \rightarrow F$ is the multiplicative character
$\det_C(x_{ij})_{N\times N} \eqdef
\det_{C \times C}(x_{ij})_{i,j \in C}$.

Hence, if $C_1, C_2, \ldots, C_r$ are the columns of $D$,
we may write
$$
Gr(D) \cong G ^r /P_D,
$$
and the bundle
$$
\LD \cong G ^r \stackrel{P_D}{\times} \det\mbox{} _D^{-1},
$$
where $P_D \eqdef P_{C_1} \times \cdots \times P_{C_r}$
and $\det_D(X_1,\ldots,X_r) \eqdef \det_{C_1}(X_1)\times \cdots
\times \det_{C_r}(X_r)$.
Under this identification,
$$
\FD \cong \mbox{ closure} \IM\left[\, \Delta G \hookrightarrow G^r
\rightarrow \Gr(D) \, \right]
$$
(c.f. Proposition \ref{bigorbit}).

Now let $\tau: G^r \rightarrow G^r$,
$\tau(g_1,\ldots,g_r) = (g_1', \ldots,g_r')$,
where $g' = \ ^t g^{-1}$, the inverse transpose of a matrix $g\in G$.
Then $\tau(P_D) = P_{D'}$, and $\tau$ induces a map
$$
\tau: \Gr(D) \rightarrow \Gr(D'),
$$
as well as a map of line bundles
$$
\begin{array}{rccc}
\tau:& \LD & \rightarrow & {\cal L}_{D'} \\
     &  \| & & \| \\
     & G^r \stackrel{P_D}{\times} \det\mbox{}_D^{-1}
	      & &
       G^r \stackrel{P_{D'}}{\times} \det\mbox{}_{D'}^{-1} \\
&(g_1, \ldots, g_r, \alpha) & \mapsto &
(g_1', \ldots, g_r', \det(g_1',\ldots,g_r') \alpha)  .
\end{array}
$$
This map is not $G$-equivariant.  Rather, if we have a
section of $\LD$,  $f:G^r \rightarrow F$ (with $f(gp) = \det_D(p) f(g)$
for $p \in P_D$), then for $g_0 \in G$,
we have
$\tau(g_0 f ) = g_0' \det(g_0')^r \tau(f)$
(a section of ${\cal L}_{D'}$).

Since $W_D$ is the restriction of such functions
$f$ to $\Delta G \subset G^r$, and $\tau(\Delta G) \subset \Delta G$,
we have an induced map
$$
\tau: W_D \rightarrow W_{D'}
$$
(an isomorphism of $F$ vector spaces),
satisfying $\tau(g_0 w) = g_0' \det(g_0')^r \tau(w)$
for $g_0 \in G$, $w\in W_D$.
This is the map required in (a),
and now (b), (c) follow trivially.
$\bullet$

\section{Resolution of singularities}
\label{resolution}

In this section,
we define the class of northwest direction diagrams,
which includes
(up to a permutation of rows and columns)
the skew, inversion, Rothe,
and column-convex diagrams.
We construct an explicit resolution of singularities
of the associated configuration varieties
by means of ``blowup diagrams''.  We also
find defining equations for these varieties.
One should note that the resolutions constructed are not
necessarily geometric blowups, and can sometimes be
small resolutions, as in Example 8 above.

We shall, as usual, think of a diagram $D$ either as
a subset of $\NN \times \NN$, or as a list $(C_1, C_2, \ldots, C_r)$
of columns $C_j \subset \NN$.
We shall examine only configuration varieties, as
opposed to line bundles on them, so we shall assume that
the columns are without multiplicity: $C_j \neq C_{j'}$
for $j \neq j'$.

\subsection{Northwest and lexicographic diagrams}

A diagram $D$ is {\em northwest} if it possesses the following
property:
$$
(i_1,j_1),\ (i_2, j_2) \in D
\Rightarrow (\min(i_1,i_2),\min(j_1,j_2)) \in D.
$$
Given two subsets $C = \{i_1 < i_2 < \ldots < i_l \}, \ \
C' = \{i'_1 < i'_2 < \ldots i'_{l'}\} \subset \NN$,
we say $C$ is {\em lexicographically less than } $C'$
\ ($C \llx C'$) if
$$
l < l\, ' \mbox{ and } i_1 = i'_1,\ \ldots ,\ i_l = i'_{l'},
$$
$$
\mbox{or } \exists\, m : \ i_1 = i'_1,\ \ldots ,\ i_{m-1} = i'_{m-1},\
i_m < i'_m.
$$
In the first case, we say $C$ is an {\em initial subset} of
$C'$ \ ($C \linit C'$).

A diagram $D = (C_1, C_2, \ldots)$ is {\em lexicographic}
if $C_1 \llx C_2 \llx \cdots$.
Note that any diagram can be made lexicographic
by rearranging the order of columns.

\begin{lem}
If $D$ is northwest, then the lexicographic rearrangement of $D$ is
also northwest.
\end{lem}

\noindent
{\bf Proof.} (a) I claim that if $j < j'$, then either
$C_j \llx C_{j'}$, or $C_{j} \ginit C_{j'}$.
Let $C_j = \{i_1 < i_2 < \ldots \}$,
$ C_{j'} = \{i'_1 < i'_2 < \ldots \}$.
We have assumed $C_j \neq C_{j'}$.
Thus $C_j \llx C_{j'}$ or
$C_j \glx C_{j'}$.
In the second case, $C_j \ginit C_{j'}$ or
there is an $r$ such that $i_1 = i'_1, \ldots  i_{r-1} = i'_{r-1},
i_r < i'_r$.
By the northwest property, this last case would mean
$i'_r \in C_j$, with $i_{r-1} = i'_{r-1} < i'_r < i_r$.
But this contradicts the definition of $C_j$.
Thus the only possibilities are those of the claim. \\
(b) It follows immediately from (a) that
if $C_1 \llx C_2 \llx \cdots \llx C_{s-1} \glx C_s$,
then there is a $t < s$ with
$C_{t-1} \llx C_s$, $C_s \linit C_t$,
$C_s \linit C_{t+1}$, \ldots,
$C_s \linit C_{s-1}$.  \\
(c)  From (b), we see that to rearrange the columns lexicographically
requires only the following operation:
we start with $C_1, C_2, \ldots$, and when we encounter the
first column $C_s$ which violates lexicographic order, we move it
as far left as possible, passing over those columns $C_i$ with
$C_s \linit C_i$.
This operation does not destroy the northwest
property, as we can easily
check on boxes from each pair of
columns in the new diagram.
By repeating this operation, we get the lexicographic rearrangement,
which is thus northwest.

\subsection{Blowup diagrams}

The  combinatorial lemmas of this section will be used to
establish geometric properties of configuration varieties.

Given a northwest diagram $D$ and two of its columns
$C, C' \subset \NN$, the {\em intersection blowup diagram}
$\Dh_{C,C'}$ is the diagram with the same columns as $D$
except that the new column
$C \cap C'$ is inserted in the proper lexicographic position
(provided $C \cap C' \neq C, C'$).

\begin{lem}
\label{intlem}
Suppose $D$ is lexicographic and northwest, and $C \llx C'$
are two of its columns. Then:
(a)  $C \cap C' \linit C'$, and
(b) if $C \subset C'$, then $C \linit C'$.
\end{lem}

\noindent
{\bf Proof.} (a) If $i \in C_j \cap C_{j'}$ and $i > i' \in C_{j'}$,
then $i' \in C_j$ by the northwest property.  Similarly for (b).
$\bullet$

\begin{lem}
If $D$ is lexicographic and northwest, then $\Dh_{C,C'}$ is also
lexicographic and northwest.
\end{lem}

\noindent
{\bf Proof.} If $C = C_j, C' = C_{j'}$ with $j < j'$, and
we insert the column $C \cap C' \linit C'$ immediately before
$C'$, then we easily check that the resulting diagram is again northwest.
Hence $\Dh_{C,C'}$, which is the lexicographic rearrangement of this,
is also northwest by a previous lemma.
$\bullet$

Consider the columns $C_1, C_2, \ldots \subset \NN$ of a northwest
diagram $D$, and take the smallest collection
$\{ \Ch_1 \llx \Ch_2 \llx \cdots \}$ of subsets of $\NN$
which contains the $C_i$ and is closed under taking intersections.
Then we define a new diagram
$\Dh = (\Ch_1, \Ch_2, \ldots )$
which we call
the {\em maximal intersection blowup diagram} of $D$.
Clearly $\Dh\, ^{\widehat{}} = \Dh$.
Repeated application of the above lemma shows that
if $D$ is lexicographic and northwest, then so is $\Dh$.

\vspace{1em}

\noindent
{\bf Examples.} For one of the (non-northwest)
diagrams considered previously,
we have:
$$
D_4 = \begin{array}{ccc}
\Box & \Box &      \\
\Box &      & \Box \\
     & \Box & \Box
\end{array}  \ \ \ \ \ \
\Dh_4 = \begin{array}{cccccc}
\Box & \Box & \Box &      &      &      \\
     & \Box &      & \Box & \Box &      \\
     &      & \Box &      & \Box & \Box
\end{array}
$$

For the diagrams $D_7$ and $D_8$ in the previous examples,
$D_8 = \Dh_7$.
$\bullet$

Consider the columns $C \subset \NN$
of a diagram $D$ as a partially ordered
set under $\subset$, ordinary inclusion.
Given two distinct columns $C$, $C'$, we say
$C'$ {\em minimally covers} $C$
(or simply $C'$ {\em covers} $C$)
if $C \subset C'$ and there is no column
of $D$ strictly included between $C$ and $C'$.

\begin{lem}
\label{maxmin}
Let $D$ be a lexicographic northwest diagram,
and $C_L$ be the last column of $D$.  Then: \\
(a) there is a column $C_l \neq C_L$ such that
$$
( \bigcup_{C \neq C_L} C)
 \cap C_L = C_l \cap C_L ;
$$
(b) if $\Dh = D$, then $C_L$ covers
at most one other column $C_l$
and is covered by
at most one other column $C_u$.
\end{lem}

\noindent
{\bf Proof.}
(a) Now, by Lemma \ref{intlem}, $C \cap C_L \linit C_L$
for any column $C$.
Hence the sets $C \cap C_L$ for $C \neq C_L$ are linearly
ordered under inclusion, and there is a largest one
$C_l \cap C_L$.
Thus
$$
( \bigcup_{C \neq C_L} C) \cap C_L
= \bigcup_{C \neq C_L} (C \cap C_L)
= C_l \cap C_L .
$$
(b) By Lemma \ref{intlem},
the columns with $C \subset C_L$ satisfy $C \linit C_L$ and are
linearly ordered, so there is at most one maximal $C_u$.

Now suppose $C_u, C_u' \llx C_L$ are columns of $D$ both
covering $C_L$.  Then again by
Lemma \ref{intlem}, we have $C_u \cap C_u' \leqlx C_u$
or $\leqlx C_u'$,
so that $C_u \cap C_u' \neq C_L$.
But
$C_u \cap C_u'$ is between $C_L$ and $C_u$,
and between $C_L$ and $C_u'$.
Hence
$C_u = C_u \cap C_u' =  C_u'$.
$\bullet$

\subsection{Blowup varieties}

Let $D = (C_1, C_2, \ldots)$ be a lexicographic northwest
diagram,
and $\Dh = (\Ch_1, \Ch_2, \ldots)$ be its maximal
intersection blowup.
Recall that $\Dh$ is obtained by adding certain columns to
$D$, so there is a natural projection map
\begin{eqnarray*}
\pr : \ \Gr(\Dh) & \rightarrow & \Gr(D) \\
	(V_{\Ch})_{\Ch \in \Dh} & \mapsto &
	(V_{C})_{C \in D},
\end{eqnarray*}
obtained by forgetting some of the linear subspaces
$V_{\Ch} \in \Gr(\Ch)$.

\begin{prop}
\label{birational}
If $D$ has $\leq N$ rows, then
$$
\pr : \ \Gr(\Dh) \rightarrow  \Gr(D)
$$
induces a birational map of algebraic varieties
$$
\pr : \ \FF_{\Dh} \rightarrow  \FD .
$$
\end{prop}

\noindent
{\bf Proof.}  Consider the dense open sets
$\Fgen_{\Dh} \subset \FF_{\Dh}$ and
$\Fgen_D \subset \FD$ of Proposition \ref{bigorbit},
consisting of subspaces in general position.
If we consider an element $g \in GL(N)$ as a sequence of
column vectors $g = (v_1,\ldots,v_n)$, and
$C = \{i_1, i_2, \ldots \} \subset \{1,\ldots,N\}$,
recall that we define
$g(C) = \Span_F(v_{i_1},v_{i_2},\ldots ) \in \Gr(C)$.
By definition, any element of $\Fgen_{\Dh}$ can
be written as
$(g(\Ch_1), g(\Ch_2), \ldots ) \in \Gr(D)$
for some $g \in GL(N)$.

Now, any column of $\Dh$ can be written as an intersection
of columns of $D$:
$\Ch = C_{j_1} \cap C_{j_2} \cap \cdots$.
Then we have
$g(\Ch) = g(C_{j_1}) \cap g(C_{j_2}) \cap \cdots$,
so the projection map
\begin{eqnarray*}
\pr : \ \Fgen_{\Dh} & \rightarrow  & \Fgen_D \\
	(g(\Ch))_{\Ch \in \Dh} & \mapsto &
	(g(C))_{C \in D}
\end{eqnarray*}
can be inverted:
$$
\begin{array}{rrcl}
\pr^{-1} : & \Fgen_D           & \rightarrow  & \Fgen_{\Dh} \\
      &  (g(C))_{C \in D} & \mapsto      & (g(\Ch) =
		   g(C_{j_1}) \cap g(C_{j_2}) \cap \cdots)_{\Ch \in \Dh}.
\end{array}
$$
Hence the map is birational on the configuration varieties as claimed.
$\bullet$

\subsection{Intersection varieties}
\label{intersection varieties}

Now, given a diagram $D$, define the {\em intersection variety}
$\ID$ of $D$ by:
$$
\ID = \{ (V_C)_{C \in D} \in \Gr(D)
	  \mid \forall C,C',\ldots \in D, \ \ \dim( V_C \cap V_{C'} \cap \cdots
)
		\geq |C \cap C' \cap \cdots|\}.
$$
Clearly $\ID$ is a projective subvariety of $\Gr(D)$,
and $\FD \subset \ID$.

If $\Dh = D$ (up to rearrangement of column order),
then the intersection conditions reduce to inclusions:
$$
\ID = \{ (V_C)_{C \in D} \in \Gr(D)
	  \mid  C \subset C' \Rightarrow V_C \subset V_{C'}
      \}.
$$

\vspace{1em}

\noindent
{\bf Example.}  For the diagram $D_4$  of Section
\ref{definitions and examples},
and $N=4$, $\II_{D_4}$ has two irreducible components,
$\FF_{D_4}$ and $\FF_{D_5}$.  That is, as before,
if we have
three lines in $\PP^3$ with non-empty  pairwise intersections,
then either they are coplanar, or they all intersect in a point.
$\bullet$

\begin{lem}
\label{union conditions}
Let $D$ be a northwest diagram, and
$\ID$ its intersection variety.
Then any configuration $(V_C)_{C \in D} \in \ID$
satisfies
$$
\dim(V_C + V_{C'} + \cdots) \leq  |C \cup C' \cup \cdots|
$$
for any columns $C, C',\ldots$ of $D$.
\end{lem}

\noindent
{\bf Proof.}
Without loss of generality, assume $D$ is lexicographic.
We use induction on the number of columns in $D$.
Now any list $C, C', \ldots$ of columns of $D$
also constitutes a lexicographic northwest diagram,
so to carry through the induction we need only
prove the statement for {\em all} the columns
$C_1, C_2, \ldots, C_L$ of $D$.
Now, by Lemma \ref{maxmin}, there is a column $C_l \neq C_L$
such that
$( \cup_{C\neq C_L} C) \cap C_L = C_l \cap C_L$.
Then we have
\begin{eqnarray*}
\dim(\ ( \sum_{C\neq C_L} V_C ) \cap V_{C_L}\ )
& \geq & \dim( \sum_{C \neq C_L} (V_C \cap V_{C_L})\  ) \\
& \geq & \dim(\, V_{C_l} \cap V_{C_L} ) \\
& \geq & |\, C_l \cap C_L| \ \ \ \mbox{  since } (V_C) \in \ID  \\
& = &  | \ ( \bigcup_{C\neq C_L} C\ ) \cap C_L \ | .
\end{eqnarray*}
Thus we may write
\begin{eqnarray*}
\dim( \sum_{C\in D} V_C ) & = &
\dim( \sum_{C \neq C_L} V_C \, ) + \dim(V_{C_L})
-\dim(\, ( \sum_{C\neq C_L} V_C \, ) \cap V_{C_L} \ ) \\
& \leq & |\! \bigcup_{C\neq C_L} C \, | + |\, C_L|
- |\, ( \bigcup_{C\neq C_L} C \ ) \cap C_L |
\ \ \ \mbox{  by induction} \\
& = & |\bigcup_{C\in D} C \ |
\ \ \ \bullet
\end{eqnarray*}

\begin{lem}
If $D$ is a northwest diagram with
$\leq N$ rows and $\Dh = D$
(up to rearrangement of column order),
then $\FD$ is an irreducible component of $\ID$.
\end{lem}

\noindent
{\bf Proof.}
Recall that $\FD$ is always irreducible.  Thus it suffices
to show that $\Fgen_D$ is an open subset of $\ID$.

Consider the set $\ID^{\mbox{\rm \tiny gen}}$ of configurations
$(V_C)_{C \in D}$ satisfying, for every list $C, C',\ldots$ of
columns in $D$,
$$
	  \dim (V_C + V_{C'} + \cdots)
	   = |C \cup C' \cup \cdots |
$$
and
$$
	  \dim (V_C \cap V_{C'} \cup \cdots)
	   = |C \cap C' \cap \cdots |   .
$$
This is an open subset of $\ID$ by the previous lemma.

I claim that $\Fgen_D = \ID^{\mbox{\rm \tiny gen}}$.
To see this equality,
let $(V_C)_{C \in D} \in \ID$ satisfy
the above rank conditions,
and we will find a basis $g = (v_1,\ldots,v_N)$ of $V=F^N$
such that $V_C = g(C)$ for all $C$. (C.f. the proof of
Proposition \ref{bigorbit}.)

As before, we consider the columns as a poset under ordinary inclusion.
We begin by choosing mutually independent bases for those
$V_C$ where $C$ is a minimal element of the poset.
This is possible because
$
\dim \Span(V_C \mid C \mbox{ minimal})
	   = \sum_{C \ \mbox{\tiny  minml}} |C|.
$

Now we consider the $V_C$ where $C$ covers a minimal column.
We start with the basis vectors already chosen, and add enough
vectors, all mutually independent, to span each space.
Again, the dimension conditions ensure there will be no conflict
in choosing independent vectors, since the $V_C$ can have no
intersections with each other except those due to the
intersections of columns.
The condition $\Dh = D$ ensures that all these intersections
are (previously considered) columns.

We continue in this way for the higher layers of the poset.
We will not run out
of independent basis vectors because
all the columns of $D$
are contained in $\{1,\ldots,N\}$.
$\bullet$

\subsection{Smoothness and equations defining varieties}

\begin{prop}
\label{smooth}
Let $D$ be a northwest diagram
with $\leq N$ rows and $\Dh = D$
(up to rearrangement of column order).
Then $\FD = \ID$, and $\FD$ is a smooth variety.
\end{prop}

\noindent
{\bf Proof.}
(a) Let $C_L$ be the last column of $D$, and
let $D'$ be $D$ without the last column.
By lemma \ref{maxmin}, $C_L$ is covered by at most one
other column $C_u$, and covers at most one other column
$C_l$.
If these columns do not exist, take
$C_l = \emptyset$, $C_u = \{1,\ldots,N\}$. \\
(b)  Now I claim that there is a fiber bundle
$$
\begin{array}{ccc}
\Gr(C_l, C_L, C_u) & \rightarrow & Z \\
		       &             & \downarrow \\
		       &             & \Gr(D') \ \,
\end{array}
$$
where $\Gr(C_l, C_L, C_u)$ denotes the Grassmannian of
$|C_L|$-dimensional linear spaces which contain a
fixed $|C_l|$-dimensional space and are
contained in a fixed
$|C_u|$-dimensional space;
and
$$
Z = \{ (\, (V_{C'})_{C'}, V_L ) \in \Gr(D') \times \Gr(C_L) \mid
V_{C_l} \subset V_L \subset V_{C_u} \}.
$$
This is clear. See also~\cite{BD1}.\\
(c)  Note that $\ID = ( \II_{D'} \times \Gr(C_L) ) \cap Z$.
This is because of the uniqueness of $C_l$ and $C_u$.
Thus the above fiber bundle restricts to
\begin{eqnarray*}
\Gr(C_l, C_L, C_u) & \rightarrow & \ID \\
		       &             & \downarrow \\
		       &             & \II_{D'} \ \ ,
\end{eqnarray*}
which is thus also a fiber bundle. \\
(d)  Now apply the above construction repeatedly,
dropping columns of $D$ from the end.
Finally we obtain $\ID$ as an
iterated fiber bundle whose fibers at each step are smooth
and connected
(in fact they are Grassmannians).
In particular, $\ID$ is smooth
and connected. \\
(e)  Since $\ID$ is a smooth, connected, projective algebraic
variety, it must be irreducible.  But by a previous lemma,
$\FD$ is an irreducible component of $\ID$.  Therefore
$\FD = \ID$, a smooth variety.
$\bullet$

\begin{prop}
\label{conn fibers}
Let $D$ be a northwest diagram with
$\leq N$ rows.  Then $\FD = \ID$, and the
birational projection map
$\FF_{\Dh} \rightarrow \FD$ has connected
fibers.
\end{prop}

\vspace{1em}

\noindent
{\bf Proof.}
(a)  I claim the following:
if $\Ch$ is a column of $\Dh$ such that
for all $C \in \Dh$ with $C \subneq \Ch$ we have
$C \in D$, then the projection map
$\II_{D \cup \Ch} \rightarrow \ID$ is onto, with
connected fibers.

Suppose $(V_C)_{C\in D}$ is a configuration in $\ID$.
Let
$$
V_u = \bigcap_{C \in D \atop C \supset \Ch} V_C
\ \ \ \
\mbox{ and } \ \ \ \
V_l = \sum_{C \in D \atop C \subset \Ch} V_C.
$$
Then $\dim(V_u) \geq |\Ch|$ since $(V_C) \in \ID$,
and $\dim(V_l) \leq |\Ch|$ by Lemma \ref{union conditions}.
Clearly $V_l \subset V_u$.
Now choose an arbitrary $V_{\Ch}$ between $V_l$ and $V_u$
with $\dim(V_{\Ch}) = |\Ch|$.
Then for any list of columns $C, C', \ldots \in D$,
we have either: \\
(i) $\Ch \cap C \cap C'\cdots = \Ch$, and
$$
V_{\Ch \cap C \cap C'\cdots} \ = \ V_{\Ch} \ = \ V_{\Ch} \cap V_u
\ \subset \ V_{\Ch} \cap V_C \cap V_{C'} \cap \cdots ;
$$
or (ii) $\Ch \cap C \cap C'\cdots \subneq \Ch$, so that
$\Ch \cap C \cap C'\cdots \in D$ by hypothesis, and
$$
V_{\Ch \cap C \cap C'\cdots} \ \subset \ V_l \cap V_C \cap V_{C'} \cdots
\ \subset \ V_{\Ch} \cap V_C \cap V_{C'} \cap \cdots .
$$

In either case $(V_C)_{C \in D \cup \Ch} \in \II_{D \cup \Ch}$.
Thus $\II_{D \cup \Ch} \rightarrow \ID$ is onto, and the
fibers are the Grassmannians $\Gr(V_l, |C|, V_u)$. \\
(b) We now see that $\II_{\Dh} \rightarrow \ID$ is onto
(with connected fibers)
by repeated application of (a), starting with
$\Ch$ minimal in the poset of columns of $\Dh$ and proceeding
upward. \\
(c)  By the previous proposition, the projection map takes
$\II_{\Dh} = \FF_{\Dh} \rightarrow \FD$.  But $\II_{\Dh} \rightarrow \ID$
is onto, so $\FD = \ID$, and we are done.
$ \bullet $

The above proposition shows that for northwest
diagrams, $\FD$ is defined by the rank conditions of $\ID$.
In general, we state the

\begin{conj}
For an arbitrary diagram $D$,
$\FD$ is the set of configurations satisfying
\begin{eqnarray*}
\dim(V_C + V_{C'} + \cdots) & \leq & |\ C \cup C' \cup \cdots| \\
\dim( V_C \cap V_{C'} \cap \cdots )
	 & \geq & |\ C \cap C' \cap \cdots|
\end{eqnarray*}
for every list $C,C',\ldots $ of columns of $D$.
Equivalently(?), the variety defined by these equations is irreducible.
\end{conj}

\section{Cohomology of line bundles}
\label{cohomology}

Using the technique of Frobenius splitting,
we show certain surjectivity and vanishing theorems for
line bundles on configuration varieties.
In particular, we show that for any
northwest diagram $D$, $\FD$ is normal, and
projectively normal with respect to $\LD$
(so that global sections of $L_D$ on $\FD$ extend to $\Gr(D)$);
and $\FD$ has rational singularities.
The material of section \ref{frobenius} was shown to me by Wilberd
van der Kallen.

\subsection{Frobenius splittings of flag varieties}
\label{frobenius}

The technique of Frobenius splitting,
introduced by V.B. Mehta,  S. Ramanan, and A. Ramanathan
{}~\cite{MR},~\cite{RR},~\cite{R1},~\cite{R2},
is a method for proving
certain surjectivity and vanishing results.

Given two algebraic varieties $Y \subset
 X$ defined
over an algebraically closed
field $F$ of characteristic $p > 0$,
with $Y$ a closed subvariety of $X$,
we say that the pair $Y \subset X$ is
{\em compatibly Frobenius split} if: \\
(i) the $p^{th}$ power map $F: \OO_X \rightarrow F_* \OO_X$
has a splitting, i.e. an $\OO_X$-module morphism
$\phi: F_* \OO_X \rightarrow \OO_X$ such that $\phi F$ is the
identity; and \\
(ii) we have $\phi(F_* I) = I$, where
$I$ is the ideal sheaf of $Y$.

Mehta and Ramanathan prove the following

\begin{thm}
\label{vanishing}
Let $X$ be a projective variety, $Y$ a closed subvariety,
and $L$ an ample line bundle on $X$.
If $Y \subset X$ is compatibly split, then
$H^i(Y,L) = 0$ for all $i >0$,
and the restriction map
$H^0(X,L) \rightarrow H^0(Y,L)$ is surjective.

Furthermore, if $Y$ and $X$ are defined and projective over $\Z$
(and hence over any field), and they are compatibly split
over any field of positive characteristic,
then the above vanishing and surjectivity
statements also hold for all fields of characteristic zero.
$\bullet$
\end{thm}

Our aim is to show that, for $D$ a northwest diagram,
$\FD \subset \Gr(D)$ is compatibly split.
The above theorem and Theorem \ref{weyl mod is image} will then imply that
$S_D^* \cong H^0(\FD,\LD) = \sum_i (-1)^i H^i(\FD,\LD)$,
the Euler characteristic of $\LD$.

We will also need the following result
of Mehta and V. Srinivas~\cite{MS}:

\begin{prop}
\label{normality}
Let $Y$ be a projective variety which is Frobenius split,
and suppose there exists
a smooth irreducible projective variety $Z$
which is mapped onto $Y$ by an algebraic map
with connected fibers.
Then $Y$ is normal.

Furthermore, if $Y$ is defined over $\Z$, and is normal
over any field of positive characteristic,
then $Y$ is also normal over all fields of characteristic zero.
$\bullet$
\end{prop}

\begin{prop}
\label{pushforward}
Let $f: Z\rightarrow X$ be a separable morphism
with connected fibers,
where $X$ and $Z$ are projective varieties
and  $X$ smooth.
If $Y \subset Z$ is compatibly split, then so is $f(Y) \subset X$.
$\bullet$
\end{prop}

We will show our varieties are split by using the above proposition
to push forward a known
splitting due to Ramanathan~\cite{R2} and O. Mathieu~\cite{M}.

For an integer $n$ and $n$ permutations $w, w', \ldots$,
define
$$
X_n = \underbrace{G \tB G \tB \cdots \tB G}_{n \
\mbox{\rm \tiny  factors}  }/B,
$$
and the twisted multiple Schubert variety
$$
Y_{w, w' \ldots}
= \overline{B w B} \tB \overline{B w' B} \tB \cdots \subset X_n
$$
Note that we have an isomorphism
$$
\begin{array}{ccc}
X_n  & \rightarrow & (G/B)^n \\
(g, g', g'', \ldots) & \mapsto & (g, g g', g g' g'',\ldots).
\end{array}
$$

\begin{prop}{(Ramanathan-Mathieu)}
\label{Ram-Math}
Let $G$ be a reductive algebraic group over a
field of positive characteristic with Weyl group $W$ and
Borel subgroup $B$,
and let $w_0, w_1, \ldots w_r \in W$.
Then $Y_{w_0, w_1, \ldots} \subset X_{r+1}$ is compatibly split.
$\bullet$
\end{prop}

Now, for Weyl group elements $u_1, \ldots, u_r$,
define a variety
$\FF_{u_1,\ldots,u_r} \subset (G/B)^r$
by
$$
\FF_{u_1,\ldots,u_r} = \overline{G \cdot ( u_1 B,\ldots, u_r B)}.
$$

\begin{prop}{(van der Kallen)}
\label{splitting}
Let $w_1, \ldots, w_r$ be Weyl group elements,
and define $u_1 = w_1, u_2 = w_1 w_2, \ldots u_r = w_1 \cdots w_r$.
Suppose $w_1, \ldots w_r$ satisfy
$\ell(w_1 w_2 \cdots w_r)
= \ell(w_1) + \ell(w_2) + \cdots + \ell(w_r)$,
or equivalently
$\ell(u_j) = \ell(u_{j-1}) + \ell(u_{j-1}^{-1} u_j)$ for all $j$.
Then the pair $\FF_{u_1,\ldots,u_r} \subset (G/B)^r$
is compatibly split.
\end{prop}

\noindent
{\bf Proof.}
Define
$$
\begin{array}{rccc}
f: & X_{r+1} & \rightarrow & (G/B)^r \\
   & (g_0, g_1, \ldots, g_r) & \mapsto &
		(g_0 g_1,\, g_0 g_1 g_2,\, \ldots,\, g_0 g_1 \ldots g_r).
\end{array}
$$
We will examine the image under this map
of
$$
Y \eqdef Y_{w_0, w_1, \ldots w_r}
= G \tB \overline{B w_1 B} \tB \cdots
\tB \overline{B w_r B} \subset X_{r+1} ,
$$
where $w_0$ is the longest element of the Weyl group.

It is well known that, under the given hypotheses,
we have
$(B w_1 B) \cdots (B w_r B) = B w_1 \cdots w_r B$,
and that the multiplication map
$$
B w_1 B \tB \cdots \tB B w_r B \rightarrow B w_1 \cdots w_r B
$$
is bijective.
Thus any element $(g, b_1 w_1 b_1',\ldots, b_r w_r b_r' B) $
(for $b_i, b_i' \in B$) can be written as
$(g , b w_1, w_2, \ldots, w_r B) = (g b, w_1, w_2, \ldots, w_r B)$
for some $b\in B$,
and
$$
f(G \tB B w_1 B \tB \cdots \tB B w_n B)
= f(G \tB w_1 \tB \ldots \tB w_r B)
= G (u_1 B, \ldots , u_r B).
$$
Hence
$f(Y)
= \FF_{u_1,\ldots,u_r}$,
since our varieties are projective.

Now, $f$ is a separable map
with connected fibers
between smooth projective varieties,
so the compatible splitting of the previous
proposition pushes forward
by Proposition \ref{pushforward}.
$\bullet$

We will need the following lemmas to show that
our configuration varieties have rational singularities.

\begin{lem}{(Kempf~\cite{K})}
\label{Kempf}
Suppose $f : Z \rightarrow X$ is a separable
morphism with generically connected fibers
between projective algebraic varieties $Z$ and $X$,
with $X$ normal.  Let $L$ be an ample line bundle on $X$,
and suppose that $H^i(Z, f^* L^{\otimes n}) = 0$
for all $i > 0$ and all $n >> 0$.

Then $R^i f_* \OO_Z = 0 $ for all $i>0$.
$\bullet$
\end{lem}

Resuming the notation of Prop \ref{Ram-Math},
let $w_1, \ldots w_n$ be arbitrary Weyl group elements,
and let $\lam_1, \ldots \lam_n$ be arbitrary weights of $G$.
Let $X_n$ be as before, and define the line bundle
$\LL_{\lam_1, \ldots \lam_n}$ on $X_n$
and on $Y_{w_1, \ldots, w_n} \subset X_n$
as the quotient of
$G^n \times F$
by the $B^n$-action
$$
(b_1, \ldots, b_n)\cdot (g_1, g_2,  \ldots, g_n, a) \eqdef
(g_1 b_1, b_1^{-1} g_2 b_2, \ldots , b_{n-1}^{-1} g_n b_n,\,
\lam_1(b_1) \cdots \lam_n(b_n) a) .
$$
Note that under the identification $X_n \cong (G/B)^n$,
$\LL_{\lam_1, \ldots, \lam_n}$ is isomorphic
to the Borel-Weil line bundle
$G^n \stackrel{B^n}{\times} (\lam_1^{-1},\ldots \lam_n^{-1})$.

\begin{lem}
\label{nonample}
Assume $\lam_1, \ldots \lam_n$ are dominant weights
(possibly on the wall of the Weyl chamber).
Then $H^i(Y_{w_1, \ldots, w_n},\, \LL_{\lam_1, \ldots, \lam_n}) = 0$
for all $i > 0$.
\end{lem}

\noindent
{\bf Proof} (van der Kallen).
Note that $\LL_{\lam_1,\ldots,\lam_n}$ is effective, but not
necessarily ample, so we cannot deduce the conclusion
directly from Theorem \ref{vanishing}.

Recall the following facts from $B$-module theory
{}~\cite{P},~\cite{vdK}: \\
(a) An excellent filtration of a $B$-module is one
whose quotients are isomorphic to Demazure modules
$H^0(\overline{BwB}, \LL_{\lam})$, for Weyl group elements $w$
and dominant weights $\lam$. \\
(b) If $M$ has an excellent filtration,
and $\EE(M) \eqdef G \tB M$ is the corresponding
vector bundle on $G/B$, then
$H^i(G/B, \EE(M)) = 0 $ for all $i >0$,
and $H^0(G/B, \EE(M))$ has an excellent filtration. \\
(c) Polo's Theorem:  If $M$ has an excellent filtration,
then so does $(\lam^{-1}) \otimes M$ for any
dominant weight $\lam$.

Now consider the fiber bundle
$$
\begin{array}{ccc}
Y_{w_2,\ldots, w_n} & \rightarrow & Y_{w_1, w_2, \ldots, w_n} \\
& & \downarrow \\
& & \overline{Bw_1B}
\end{array}
$$
which leads to the spectral sequence
$$
H^i(\  \overline{Bw_1B}, \
\EE( \, (\lam_1^{-1}) \otimes
     H^j(Y_{w_2,\ldots , w_n}, \LL_{\lam_2,\ldots, \lam_n})\, )\  )
\Rightarrow
H^{i+j}(Y_{w_1, w_2, \ldots, w_n}, \LL_{\lam_1, \lam_2,\ldots, \lam_n}) .
$$
By induction, assume that
$H^j(Y_{w_2, \ldots, w_n}, \LL_{\lam_2,\ldots,\lam_n}) = 0$ for $j>0$,
and that
$H^0(Y_{w_2, \ldots, w_n}, \LL_{\lam_2,\ldots, \lam_n})$
has an excellent filtration.
Then applying (b) and (c), we find
$$
H^i(Y_{w_1, w_2, \ldots, w_n}, \LL_{\lam_1, \lam_2,\ldots, \lam_n}) =
H^i(\ \overline{Bw_1B},
       \, \EE( \, (\lam_1^{-1}) \otimes
	     H^0(Y_{w_2,\ldots , w_n}, \LL_{\lam_2,\ldots, \lam_n})\, ) \ )
 = 0
$$
for $i>0$,
and that
$H^0(Y_{w_1, w_2, \ldots, w_n}, \LL_{\lam_1, \lam_2,\ldots, \lam_n})$
has an excellent filtration.
$\bullet$

\subsection{Frobenius splitting of Grassmannians}

We would now like to push forward the
Frobenius splittings found
above for flag varieties
to get splittings of configuration varieties.
For this we need a combinatorial prerequisite.

Given a diagram $D = (C_1, C_2,\ldots, C_r)$
with $\leq N$ rows,
consider a sequence of permutations
(Weyl group elements)
$u_1, u_2, \ldots \in \Sigma_N$
such that, for all $j$: \\
($\alpha$) $\ell(u_j) = \ell(u_{j-1})+\ell(u_{j-1}^{-1} u_{j}) $, and \\
($\beta$) $u_j(\ \{1,2,\ldots,|C_j| \ \}) = C_j$.
The first condition says that the sequence is
increasing in the weak order on the Weyl group.
In the next section, we will give an algorithm
which produces such a sequence for
any northwest diagram, so that the following
theorem will apply:

\begin{prop}
\label{config split}
If $D$ a diagram which admits a sequence
of permutations $u_1, u_2, \ldots$ satisfying
($\alpha$) and ($\beta$) above,
then the pair $\FD \subset \Gr(D)$ is compatibly split for any
field $F$ of positive characteristic.

\mbox \\
Hence over an algebraically closed field $F$ of
arbitrary characteristic, \\
(a) the cohomology groups $H^i(\FD, \LD) = 0$ for $i > 0$; \\
(b)  the restriction map $\rest_{\Delta} : H^0(\Gr(D), \LD)
\rightarrow H^0(\FD, \LD)$ is surjective;\\
(c)  $\FD$ is a normal variety.
\end{prop}

\noindent
{\bf Proof.}
By ($\beta$), the maximal parabolic
subgroups $P_C = \{ (x_{ij}) \in GL(N) \mid
x_{ij} = 0 \mbox{ if } i \not\in C , \ j \in C \}$
satisfy
$u_i B u_i^{-1} \subset P_{C_i}$.
Write
$$
\Gr(D)  = \Gr(C_1) \times \cdots \times \Gr(C_r)
\cong G/P_{C_1} \times \cdots \times G/P_{C_r},
$$
and consider the $G$-equivariant projection
$$
\begin{array}{rccc}
\phi: & (G/B)^r & \rightarrow & \Gr(D) \\
   & (g_1 B,\ldots,g_r B) & \mapsto &
	(g_1 u_1^{-1} P_{C_1},\ldots,g_r u_r^{-1} P_{C_r})
\end{array}
$$
Then we have $\phi(u_1 B, \ldots, u_r B) = (I \, P_{C_1},\ldots, I \, P_{C_r})$
and $\phi(\FF_{u_1,\ldots,u_r}) = \FD$.
Since $\phi$ is a map with connected fibers
between smooth projective varieties,
we can push forward the compatible splitting
for $\FF_{u_1,\ldots,u_r} \subset (G/B)^r$ found in the
previous section.  Applying
Theorem \ref{vanishing} and Propositions \ref{normality} and
\ref{conn fibers},
we have the assertions of the theorem.
$\bullet$

Note that (b) and (c) of the Proposition are equivalent
to the projective normality of $\FD$ with respect to $\LD$.

\begin{conj}
For any diagram $D$, and any Weyl group elements
$u_1, \ldots u_r$, the pairs $\FD \subset \Gr(D)$ and
$\FF_{u_1, \ldots u_r} \subset (G/B)^r$ are compatibly
split.
\end{conj}

In order to prove the character formula in the last
section of this paper, we will need stronger relations
between the singular configuration varieties and their
desingularizations.  In particular,
we will show that our varieties have
rational singularities.

\begin{lem}
Let $X$, $Y$ be algebraic varities with an action of an
algebraic group $G$, and $f: X\rightarrow Y$ an equivariant
morphism.  Assume that $X$ has an open dense $G$-orbit $G\cdot x_0$,
and take $y_0 = f(x_0)$, $G_0 = \mathop{\rm Stab}_G y_0$.

Then $f^{-1}(y_0) = \overline{G_0 \cdot x_0}$.
In particular, if $G_0$ is connected,
then $f^{-1}(y_0)$ is connected and irreducible.
\end{lem}

\noindent
{\bf Proof.}
For $F = \C$, this is trivial.  Take $x_1 \in f^{-1}(y_0)$, and
consider a path $x(t) \in X$ such that $x(0) = x_1$ and $x(t) \in G\cdot x_0$
for small $t>0$.  Then the path $f(x(t))$ lies in $G\cdot y_0$ for
small $t \geq 0$, and we can lift it to a path $g(t) \in G$ such that
$g(0) = \id$ and $f(x(t)) = g(t)\cdot y_0$ for small $t \geq 0$.
Then $\tilde{x}(t) \eqdef g(t)^{-1}\cdot x(t)$ satisfies
$\tilde{x}(0) = x_1$, $\tilde{x}(t) \in G_0\cdot x_0$ for small
$t > 0$.

For general $F$, T. Springer
has given the following clever argument.
Assume without loss of generality that $X$ is irreducible and
$G\cdot y_0$ is open dense in $Y$.  Since an algebraic map
is generically flat, and $G\cdot y_0$ is open, all the irreducible
components $C$ of $f^{-1}(y_0)$ have the same dimension
$\dim C = \dim X - \dim Y$.   Let $Z = \overline{G\cdot C}$
be the closure of one of these components.
Now, the restriction $f: Z \rightarrow Y$ also satisfies our
hypotheses, with $C \subset Z$ again a component of
the fiber of the restricted $f$,
so we again have $\dim C = \dim Z - \dim Y$, and
$\dim Z = \dim X$. Thus $G\cdot C$ is an open subset
of $X$, since $X$ is irreducible.

Now consider the open set $G\cdot C \cap G\cdot x_0 \subset X$.
Choose a point $z$ in this set which does not lie in any other
component $C'$ of our original $f^{-1}(y_0)$.  For any other
component $C'$, choose a similar point $z'$.  But we have
$g \cdot z \in C$,
$g'\cdot z' = g_0 g\cdot z \in C'$ for some $g, g', g_0 \in G$,
and in fact $g_0 \in G_0$.
Thus $C' = g_0 \cdot C$, and $G_0$ permutes the components transitively.
Hence, $G_0 \cdot x_0 $ has at least
as many irreducible components as the whole $f^{-1}(y_0)$,
and the lemma follows.
$\bullet$

\begin{prop}{(Inamdar-van der Kallen)}
\label{rational sing}
Suppose $D_1$, $D_2$ are diagrams admitting sequences
of permutations with ($\alpha$) and ($\beta$) as above, such that
$D_2$ is obtained by removing some of the columns of $D_1$.
Denote $\FF_1 = \FF_{D_1}$, $\FF_2 = \FF_{D_2}$, $\LL_2 = \LL_{D_2}$,
and consider the projection
$\pr : \FF_1 \rightarrow \FF_2$.\\
Then: \\
(a) $H^0(\FF_1, \pr^* \LL_2 ) =
H^0(\FF_2, \LL_2 )$, and this $G$-module has a good
filtration (one whose quotients are isomorphic to
$H^0(G/B, \LL_{\lam})$ for dominant weights $\lam$).
 \\
(b) $H^i(\FF_1, \pr^* \LL_2 ) =
H^i(\FF_2, \LL_2 ) = 0$ for all $i > 0$.
 \\
(c)  $R^i \pr_* \OO_{\FF_2} = 0$ for all $i >0$.
\\
(d)  If $F$ has characteristic zero, then
$\FD$ has regular singularities for any
northwest diagram $D$.
\end{prop}

\noindent
{\bf Proof.}
(i) Consider a sequence of permutations $w_1, w_2, \ldots w_r$
(where $r$ is the number of columns in $D_1$) such that
$u_1 = w_1, u_2 = w_1 w_2, \ldots $ satisfies
($\alpha$) and ($\beta$), and let $Y = Y_{w_0, w_1, \ldots w_r}$,
(where $w_0$ is the longest permutation).
Then we have a commutative
diagram of  surjective morphisms
$$
\begin{array}{ccc}
Y & \stackrel{\Phi_1}{\rightarrow} & \FF_1  \\
  & \stackrel{\Phi_2}{\searrow}  & \downarrow \mbox{\rm \small pr}    \\
  &                                & \FF_2
\end{array}
$$
where $\Phi_j = \phi \, \circ f$, where $\phi$ and $f$ are the
maps defined in the proofs of Propositions \ref{splitting} and
\ref{config split}
in the cases $D = D_j$.
All of these spaces have dense $G$-orbits.
Furthermore, the stabilizer of a general point in $\FD$ is an
intersection of parabolic subgroups and is connected.
Thus, by the above lemma, the fibers of $\Phi_1$ are generically
connected. \\
(ii)  Now (i) and Lemma \ref{nonample} insure that the hypotheses
of Kempf's lemma (Proposition \ref{Kempf}) are satisfied.
Thus $R^i(\Phi_1)_* \OO_Y = 0$ for $i>0$,
and by the Leray spectral sequence
we have, for all $i \geq 0$,
$$
H^i(Y, \Phi_1^* \pr^* \LL_2) =
H^i(\FF_1, (\Phi_1)_* (\Phi_1)^* \pr^* \LL_2) .
$$
(iii)  Furthermore, $\FF_1$ is normal by the previous Proposition,
and $\Phi_1$ is separable with connected fibers, so
\begin{eqnarray*}
 (\Phi_1)_* (\Phi_1)^* \pr^* \LL_2 & \cong &
[ (\Phi_1)_* (\Phi_1)^* \OO_{\FF_1} ] \otimes \pr^* \LL_2 \\
 & \cong & \pr^* \LL_2.
\end{eqnarray*}
Thus $H^i(Y, \Phi_1^* \pr^* \LL_2) = H^i(\FF_1, \pr^* \LL_2)$
for all $i \geq 0$. \\
(iv)  An exactly similar argument shows that
$H^i(Y, \Phi_2^* \LL_2) = H^i(\FF_2, \LL_2)$ for all $i \geq 0$.
But $\Phi_2^* = \Phi_1^* \pr^*$, so for all $i$,
$$
H^i(\FF_1, \pr^* \LL_2) = H^i(Y, \Phi_2^* \LL_2) =
H^i(\FF_2, \LL_2) .
$$
But we saw in Lemma \ref{nonample}
that $H^i(Y, \Phi_2^* \LL_2)$ vanishes for $i>0$,
so (b) of the present Proposition follows. \\
(iv)  We also saw in the proof of Lemma \ref{nonample}
that $H^0(Y, \Phi_2^* \LL_2)$ has an excellent filtration
as a $B$-module.  But this is equivalent to it having
a good filtration as a $G$-module, so (a) follows. \\
(v)  Now consider the spectral sequence
$$
R^i \pr_* R^j(\Phi_1)_* \OO_Y \Rightarrow R^{i+j} (\Phi_2)_* \OO_Y.
$$
For $i>0$, we have $ R^i(\Phi_1)_* \OO_Y = 0$ and
$R^i(\Phi_2)_* \OO_Y =0 $ by (ii) above.
Because of this and the normality of $\FF_1$,
we have for all $i>0$,
\begin{eqnarray*}
0 & = & R^i(\Phi_2)_* \OO_Y \\
  & = & R^i\pr_* (\Phi_1)_* \OO_Y \\
  & = & R^i \pr_* \OO_{F_1} .
\end{eqnarray*}
this shows (c). \\
(vi)  Now take $D_2 = D$ an arbitrary northwest diagram, and
$D_1 = \Dh$ its maximal blowup.  Then $\pr$ is a resolution of
singularities by Proposition \ref{smooth}.
Assume, as we will show in the next section,
that $D_1$ admits a sequence of permutations as required.
Then (c) holds, and this is precisely the definition of
rational singularities in characteristic zero, so we have (d).
$\bullet$

\subsection{Monotone sequences of permutations}

Let $D = (C_1, C_2,\ldots, C_r)$
be a northwest diagram
with $\leq N$ rows.
In this section, we will construct by a recursive
algorithm  a sequence of permutations
$u_1, u_2, \ldots \in \Sigma_N$
satisfying the conditions
of the previous section: for all $j$, \\
($\alpha$) $\ell(u_j) = \ell(u_{j-1})+\ell(u_{j-1}^{-1} u_{j}) $, and \\
($\beta$) $u_j(\ \{1,2,\ldots,|C_j| \ \}) = C_j$.

For each column $C$ of
$D$, define
the integer
$$
\gap_   N(C) = \left\{
\begin{array}{l}
 \max \{ i\mid i\not\in C,\ \exists i'\in C : i < i' \}
\mbox{, if this set is } \neq \emptyset \\
 N \mbox{, if the above set is empty.}
\end{array}
\right.
$$
Since $D$ is northwest, there is an integer $J_N \geq 1$ such that
$$
N = \gap_   N(C_1) = \cdots = \gap_   N(C_{J_N-1}) >
\gap_   N(C_{J_N}) = \cdots = \gap_  N(C_r).
$$

Now define the {\em derived diagram} $D'$ of $D$
as follows.  Given a column $C$ of $D$,
there is a  corresponding column $C'$ of $D'$:
$$
C' = \{i\mid i \in C, \ i < \gap_  N(C) \}
	     \cup \{ i-1 \mid i\in C,\ i > \gap_   N(C) \}.
$$
That is, we take $C$ and push
all squares below the $\gap_   N(C)$-th row upward
by one place.

\begin{lem}
If $D$ is northwest with $\leq N$ rows, then $D'$ is  northwest
with $\leq N-1$ rows.
\end{lem}

\noindent
{\bf Proof.}
The only doubtful case in checking the northwest property is
that of two squares $(i_1,j_1)$ and $(i_2,j_2)$ in $D'$
with $j_1 < J_N \leq j_2$ and $i_1 > i_2$.
Since $j_1 < J_N$, we have $C_{j_1} = \{1,2,\ldots,i_1-1,i_1,\ldots\}$,
so that $i_2 \in C_{j_1}$ and $i_2 \in C_{j_1}'$.
Hence $(i_2,j_1) \in D'$
as required.
$\bullet$

Now, consider the following elements of $\Sigma_N$:
$$
\kappa_n^{(N)}(i) = \left\{ \begin{array}{ll}
i & \mbox{ if } i<n \\
i+1 & \mbox{ if } n\leq i <N \\
n & \mbox{ if } i = N
\end{array}
\right.
$$
Then $\kappa_1^{(N)},\ldots,\kappa_N^{(N)}$ are
minimum length coset representatives
of the quotient $\Sigma_N/ \Sigma_{N-1}$,
and for any permutation $\pi \in \Sigma_{N-1}$,
we have $\ell(\kappa_n \pi) = \ell(\kappa_n) + \ell(\pi)$.

Now, starting with $D$, a northwest diagram with $\leq N$ rows,
we can define a sequence of derived diagrams
$D = D^{(N)}, D^{(N-1)}, \ldots, D^{(1)}$,
where $D^{(i)} = (D^{(i+1)})'$ is a northwest diagram
with $\leq i$ rows.
Let the columns of $D^{(i)}$ be $C_1^{(i)},\ldots,C_r^{(i)}$,
and define $\gap(i,j) = \gap_  i(C_j^{(i)})$.
For each $i$, we have
$$
i = \gap(i,1) = \cdots = \gap(i,J_i-1)
> \gap(i,J_i) = \cdots = \gap(i,r).
$$
Then either $\kappa_{\gap(i,j)}^{(i)}(\{1, 2, \ldots,i-1\})
\supset C_j^{(i)}$,
or $C_j^{(i)} = \{1, 2,\ldots, i\}$.
Notice that $J_N \leq J_{N-1} \leq \cdots$.

Finally, for each column $j = 1,\ldots, r$, define
$$
u_j = \kappa_{\gap(N,j)}^{(N)}\,
\kappa_{\gap(N-1,j)}^{(N-1)}
\cdots  \kappa_{\gap(1,j)}^{(1)}.
$$
This is a reduced decomposition of $u_j$, in the sense
that $\ell(u_j)$ is the sum of the lengths of the factors.
Since $\kappa_i^{(i)} = \id$, and $J_N \leq J_{N-1} \leq \cdots$,
each $u_j$ is an initial string of $u_{j+1}$.
Thus the $u_j$ have the desired monotonicity property ($\alpha$).

It only remains to show property ($\beta$):

\begin{lem}
For each column $C_j$ of $D$,
$u_j(\{1,2,\ldots,|C_j|\}) = C_j$.
\end{lem}

\noindent
{\bf Proof.}
For each $i$, we have a $u_j^{(i)}$ associated to $D^{(i)}$,
with $u_j^{(i+1)} = \kappa_{\gap(i+1,j)}^{(i+1)}\, u_j^{(i)}$.

For a given $i<N$, assume that
$u_j^{(i)}(\{1,2,\ldots, |C_j^{(i)}|\}) \subset C_j^{(i)}$.
Then I claim the same is true for $i+1$.

This is clear, because $C_j^{(i+1)}$ is $C_j^{(i)}$ with some
of its squares pushed down, and $\kappa_{\gap(i+1,j)}^{(i+1)}$
pushes down these squares to the proper positions.
In the case that $C_j^{(i+1)} = \{1,2,\ldots,l\}$, for some $l$,
we have
$u_j^{(i+1)} = u_j^{(i)} = \id$, and the claim is again true.

The lemma now follows by induction on $i$.
$\bullet$

\section{A Weyl character formula}

The results of the last two sections
allow us to apply the Atiyah-Bott Fixed Point Theorem
to compute the characters of the
Schur modules for northwest diagrams.
To apply this theorem, we must examine the points of
$\FD$ fixed under the action of $H$, the group of
diagonal matrices.  We must also understand the action
of $H$ on the tangent spaces at the fixed points.

\subsection{Fixed points and tangent spaces}
\label{fixed points}

The following formula is due to Atiyah and Bott~\cite{AB}
in the complex analytic case, and was extended to the
algebraic case by Nielsen~\cite{N},~\cite{Iv}.

\begin{thm}
\label{AB thm}
Let $F$ be an algebraically closed field,
and suppose the torus $H = (F^{\times})^N$
acts on a smooth projective
variety $X$ with isolated fixed points,
and acts equivariantly on a line bundle $L \rightarrow X$.
Then the character of $H$ acting on the cohomology groups of $L$
is given by:
$$
\sum_i (-1)^i \tr(h \mid H^i(X,L)) =
\sum_{p \ \mbox{\tiny  fixed}} {\tr(h \mid L|_p )
\over  \det( \id -\, h \mid T^*_p X )},
$$
where $p$ runs over the fixed points of $H$, $L|_p$ denotes the
fiber of $L$ above $p$, and $T^*_p X$ is the cotangent space.
$\bullet$
\end{thm}
We will apply the formula for $X = \FD$ a smooth configuration variety,
where $D = (C_1, C_2, \ldots)$ is a lexicographic northwest
diagram  with $\leq N$ rows and $\Dh = D$.

\vspace{1em}

\noindent
{\bf Fixed points.}
Assume for now that the columns are all distinct.
Let
$H = \{h=\diag(x_1,\ldots,x_N) \in GL(N) \}$
act on $\Gr(D)$ and $\FD$ by the restriction of the $GL(N)$ action.
Then by Proposition \ref{smooth}, we have
$\FD = \ID = \{(V_C)_{C \in D} \in \Gr(D) \mid
      C \subset C' \Rightarrow V_C \subset V_{C'}  \}$,
a smooth variety.

A point in $\FD \subset \Gr(D) = \Gr(C_1) \times \Gr(C_2) \times \cdots$
is fixed by $H$ if and only if each component is fixed.
Now, the fixed points of $H$ in $\Gr(l, F^N)$ are
 the coordinate planes
$E_{k_1,\ldots,k_l} = \Span(e_{k_1}, \ldots, e_{k_l})$, where
the $e_k$ are coordinate vectors in $F^N$
(c.f.~\cite{H}).
For instance, the fixed points in $\PP^{N-1}$ are the $N$ coordinate
lines $F e_k$.
We may describe the fixed points in $\Gr(C)$ as
$E_S = \Span(e_k \mid k\in S)$, where
$S \subset \{1,\ldots,N\}$ is any set with $|S| = |C|$.

Hence the fixed points in $\FD$ are as follows:
Take a function $t$ which assigns to any column $C$
a set $t(C)\subset \{1,\ldots,N\}$ with $|\, t(C)| = |C|$,
and $C \subset C' \Rightarrow t(C) \subset t(C')$.
We will call such a $t$ a {\em standard column tabloid}
for $D$.  Then the fixed point corresponding to $t$ is
$E_t = (\, E_{t(C)})_{C\in D}$.

\vspace{1em}

\noindent
{\bf Tangent spaces at fixed points.}
We may naturally identify the tangent space
$T_{V_0} Gr(l, F^N) = \Hom_F(V_0, F^N/V_0)$.
If $V_0$ is a fixed point (that is, a space stable under $H$),
then $h \in H$ acts
on a  tangent vector $\phi \in Hom_F(V_0,F^N/V_0)$
by $(h\cdot \phi)(v) = h (\phi( h^{-1} v))$.
For $(V_C)_{C\in D} \in \Gr(D)$, we have
$T_{(V_C)} \Gr(D) = \bigoplus_{C\in D} Hom(V_C, F^N/V_C)$.
Furthermore, if $(V_C)_{C\in D} \in \FD$, then
$$
\begin{array}{cl}
T_{(V_C)} \FD  = \{ \phi  = (\phi_C)_{C\in D} &
\! \! \in \bigoplus_{C\in D} Hom(V_C, F^N/V_C)  \mid  \\
& C \subset C' \Rightarrow \phi_{C'}|_{V_C} \equiv
\phi_C \mbox{ mod } V_{C'} \}
\end{array}
$$
(that is, the values of $\phi_{C}$ and $\phi_{C'}$ on
$V_C$ agree up to translation by elements of $V_{C'}$).
See~\cite{H}.

For a fixed point $E_t = (E_{t(C)})$, we will find a basis
for $T_{E_t}$ consisting of eigenvectors of $H$.
Now, the eigenvectors in $T_{E_t} \Gr(D) =
\bigoplus_{C\in D} Hom(V_C, F^N/V_C)$ are precisely
$\phi^{ijC_0}= (\phi^{ijC_0}_C)_{C \in D}$,
where $i,j\leq N$, $C_0$ is a fixed column of $D$,
and $\phi^{ijC_0}_C(e_l) \eqdef \delta_{C_0, C} \delta_{il} e_j$
($\delta$ being the Kronecker delta).
The eigenvalue is
$$
h \cdot \phi^{ijC_0} =  \diag(x_1, \ldots x_N)
\cdot \phi^{ijC_0} = x_i^{-1} x_j \ \phi^{ijC_0}.
$$

To obtain eigenvectors of $T_{E_t} \FD$, we must impose the
compatibility conditions.  An eigenvector $\phi$ with eigenvalue
$x_i^{-1} x_j$ must be a linear combination
$\phi = \sum_{C \in D} a_C \phi^{ijC}$ with $a_C \in F$.
By the compatibility, we have  that
$$
C\subset C',\, i\in t(C),\, j\not\in t(C')\ \ \Rightarrow \ \ a_C = a_{C'} .
$$
We wish to find the number $d_{ij}$
of linearly independent solutions of this condition
for $a_C$.

Given a poset with a relation $\subset$,
define its {\em connected components}
as the equivalence classes generated
by the elementary relations $x \sim y$ for $x \subset y$.
Now for a given $i,j$ consider the poset whose elements are those
columns $C$ of $D$ such that $i\in t(C)$, $j \not\in t(C)$,
with the relation of ordinary inclusion.
Then $d_{ij}$ is the number of components of this poset.

Note that the eigenvectors for all the eigenvalues
span the tangent space.
Thus
$$
\det(\, \id - h \ | \ T^*_{E_t}) =  \prod_{i\neq j}
(1 - x_i x_j^{-1})^{d_{ij}(t)}.
$$

\vspace{1em}

\noindent
{\bf Bundle fibers above fixed points.}
Finally, let us examine the line bundles $L$
on $\FD$  obtained by giving each column $C$
a multiplicity  $m(C) \geq 0$.
If $m(C) > 0$ for all columns $C$ of $D$, then $L \cong \LL_{D'}$
for the diagram $D'$ with the same columns as $D$, each
repeated $m(C)$ times.  If some of the $m(C) = 0$,
then $L$ is the {\em pullback} of $\LL_{D'}$
for the diagram
$D'$ with the same columns as $D$, each taken $m(C)$ times,
where 0 times means deleting the column.
In the second case, $L$ is effective, but not ample.

It follows easily from the definition that
$$
\tr( h \mid L|_{E_t} ) = x_1^{-\wt _1(t)} \cdots x_N^{-\wt _N(t)},
$$
where
$$
\wt_ i(t) = \sum_{C \atop i \in t(C)} m(C) .
$$

Hence we obtain:
$$
\sum_i (-1)^i \tr(h \mid H^i(\FD,\LD)) =
\sum_{t} {  \prod_i x_i^{-\wt_ i(t)} \over
\prod_{i\neq j} (1-x_i x_j^{-1})^{d_{ij}(t)}  },
$$
where $t$ runs over the standard column tabloids of $D$.

\subsection{The character formula}
\label{character formula}

We summarize in combinatorial language the implications of the
previous section.

We think of a diagram $D$
as a list of columns $C_1, C_2, \ldots \subset \NN$,
possibly with repeated columns.
Given a diagram $D$, the {\em blowup diagram} $\Dh$ is the
diagram whose columns consist of all the columns of $D$ and all
possible intersections of these columns.  We will call the
columns which we add to $D$ to get $\Dh$ the {\em phantom columns}.

We may define a {\em standard column tabloid} for the diagram $\Dh$\,
with respect to $GL(N)$,
to be a filling (i.e. labeling)
of the squares of $\Dh$ by integers in $\{1,\ldots,N\}$,
such that:\\
(i)  the integers in each column are strictly increasing, and\\
(ii) if there is an inclusion  $C \subset  C'$ between two columns,
then all the numbers in the filling of $C$ also appear in the filling of $C'$.

Given a tabloid $t$ for $\Dh$, define integers $\wt_ i(t)$ to be the
number of times $i$ appears in the filling, but
{\em not counting i's which appear in the phantom columns}.
Also define integers $d_{ij}(t)$ to be the number of connected
components of the following graph:  the vertices are columns $C$ of
$\Dh$ such that $i$ appears in the filling of $C$, but $j$ does not;
the edges are $(C,C')$ such that $C \subset C'$ or $C' \subset C$.
(An empty graph has zero components.)

Recall that a diagram $D$ is {\em northwest} if
$ i \in C_j,\ i' \in C_{j'} \ \Rightarrow \min(i,i') \in C_{\min(j,j')} $.
The following theorem applies without change to {\em northeast} diagrams
and any other diagrams obtainable from northwest ones
by rearranging the order of the rows
and the order of the columns.
Also, we can combine it with Theorem \ref{complement thm}  to
compute the character for the complement of
a northwest diagram in an $N \times r$ rectangle.

Denote a diagonal matrix by
$h = \diag(x_1,\ldots,x_N).$

\begin{thm}
Suppose $D$ is a northwest diagram with $\leq N$ rows, and
$F$ an algebraically closed field.  Then:\\
(a) The character of the Weyl module $W_D$ (for $GL(N,F)$) is given by
$$
\Char_  {W_D}(h) = \sum_{t} {  \prod_i x_i^{-\wt_ i(t)}
\over \prod_{i\neq j} (1-x_i x_j^{-1})^{d_{ij}(t)}   },
$$
where $t$ runs over the standard tabloids for $\Dh$.\\
(b) For $F$ of characteristic zero, the character of the Schur module
$S_D$ (for $GL(N,F)$) is given by
$$
\Char_  {S_D}(h) = \sum_{t} {  \prod_i x_i^{\wt_ i(t)}
\over \prod_{i\neq j} (1-x_i^{-1} x_j)^{d_{ij}(t)}  },
$$
where $t$ runs over the standard tabloids for $\Dh$.
\end{thm}

\vspace{1em}

\noindent
{\bf Example.}
Consider the following diagram
and some of its standard tabloids for $N = 3$:
$$
D = \begin{array}{ccc}
 \Box & \Box &      \\
      & \Box & \Box
\end{array}
\ \ \ \
t_1 = \begin{array}{ccc}
 1 & 1 &   \\
   & 2 & 2
\end{array}
\ \ \ \
t_2 = \begin{array}{ccc}
 1 & 1 &   \\
   & 2 & 1
\end{array}
\ \ \ \
t_3 = \begin{array}{ccc}
 3 & 2 &   \\
   & 3 & 3
\end{array}
$$
The tabloid $t_1$ has $d_{12} = d_{13} = d_{21} = d_{23} = 1$,
$d_{31} = d_{32} = 0$, and $t_2$ has $d_{12} = 2$, $d_{13} = d_{23} = 1$,
$d_{21} = d_{31} = d_{32} = 0$.  The other standard tabloids can be
obtained from $t_1$ and $t_2$ by applying a permutation of
$\{1,2,3\}$ to the entries, and rearranging the entries in the middle
column to make them increasing.  For instance, $t_3 = \pi t_2$, where
$\pi$ is the transposition $(13)$.
Note that the standard tabloids are  standard
tableaux in the usual sense:  they are fillings
with the columns strictly increasing, and the rows non-increasing.
This is true in general when $\Dh$ is a skew diagram with no repeated
columns, though not all the standard tableaux are obtained in this way.

Applying our formula we find that
$ \Char \, S_D = s_{(3,1,0)} + s_{(2,2,0)}$,
where $s_{(\lam_1,\lam_2,\lam_3)}$ is a classical Schur function,
the character of an irreducible Schur module.
Since $D$ is a skew diagram,
we could have obtained this result
using the Littlewood-Richardson Rule.
It should be possible to prove this rule using
the present methods.
$\bullet$

\vspace{1em}

\noindent
{\bf Proof of the Theorem.}
(i) Consider the map $\pr : \FF_{\Dh} \rightarrow \FD$, and the
pullback line bundle $\pr ^* \LD$.  This is the bundle on $\FF_{\Dh}$
corresponding to giving the phantom columns $C$ of $\Dh$ multiplicity
$m_C = 0$.
Let RHS denote the right hand side of our formula in (a).
Then by the analysis of Section \ref{fixed points},
RHS is equal to the right hand
side of the Atiyah-Bott formula (Theorem \ref{AB thm}) for
$X = \FF_{\Dh}$, $L = \pr ^* \LD$.
Thus
$$
\mbox{RHS} = \Char \sum_i (-1)^i H^i(\FF_{\Dh}, \pr ^* \LD) .
$$
(ii) By Proposition \ref{rational sing},
we have
$ H^i(\FF_{\Dh}, \pr^* \LD) = 0 $ for $i>0$,
and $ H^0(\FF_{\Dh}, \pr^* \LD) = H^0(\FD, \LD)$.
Thus RHS = $\Char H^0(\FD,\LD)$. \\
(iii)  By Proposition \ref{config split},
the restriction of global sections of $\LD$
from $\Gr(D)$ to $\FD$ is  surjective,
and we have
\begin{eqnarray*}
\mbox{RHS} & = &\Char H^0(\FD, \LD) \\
& = &
\Char
\IM\left( \rest : H^0(\Gr(D),\LD) \rightarrow H^0(\FD,\LD) \right) \\
& = & \Char W_D.
\end{eqnarray*}
The last equality holds by
Proposition \ref{weyl mod is image}, and we have proved (a).
Then (b) follows because $S_D = (W_D)^*$.
$\bullet$

\subsection{Betti numbers}

In this section, we compute the betti numbers of the
smooth configuration varieties
of Section \ref{fixed points}.

\begin{prop}{(Bialynicki-Birula~\cite{B})}
Let $X$ be a smooth projective variety over an algebraically
closed field $F$, acted on by the one-dimensional torus $F^{\times}$
with isolated fixed points.  Then there is a decomposition
$$
X = \coprod_{p \  \mbox{\rm \tiny  fixed}} X_p ,
$$
where the $X_p$ are disjoint, locally closed, $H$-invariant
subvarieties, each isomorphic to an affine space
$X_p \cong \mbox{\rm \bf A}^{d^+\! (p)}$.

The dimensions $d^+(p)$ are given as follows.
Let the tangent space $T_p X \cong \bigoplus_{n \in \Z} a_n(p) F_n$, where
$a_n(p) \in \NN$ and $F_n$ is the one-dimensional representation
of $F^{\times}$ for which the group element $t \in F^{\times}$
acts as the scalar $t^n$.
Then
$$
d^+\!(p) = \sum_{n > 0} a_n(p) .
$$
$\bullet$
\end{prop}

Over $\C$, the above proposition does not quite give a CW decomposition
for $X$, since the boundaries of the cells need not lie in cells of lower
dimension.  Nevertheless, $\dim_{\R} \partial X_p \leq \dim_{\R} X_p - 2$,
and this is enough to fix the betti numbers $\beta_i = \dim_{\R} H^{i}(X,
\R)$:\ \
$\beta_{2i} = \#\{p \mid d^+(p) = i\}$, and $\beta_{2 i + 1} = 0$.

Now, in our case consider the spaces $X = \FD$ of Section \ref{fixed points},
acted on by the $N$-dimensional torus $H$.  Consider the embedding
$$
\begin{array}{cccc}
\check{\rho} : & F^{\times} & \rightarrow & H \\
 & t & \mapsto &  \diag(t^{N-1}, t^{N-2}, \ldots, t, 1),
\end{array}
$$
corresponding to the
coweight $\check{\rho} = $ half sum of the positive coroots.
Then $\check{\rho}(F^{\times})$ has the same (isolated)
fixed points as $H$, since none of the
eigenvectors of $H$ on $T_p X$ is fixed
by $\check{\rho}(F^{\times})$.
(I.e., $(\alpha, \check{\rho}) \neq 0$ for any root $\alpha$.)

Also, a given eigenvector of weight $x_i x_j^{-1}$ is
of positive $\check{\rho}(F^{\times})$ weight
exactly when $i < j$.
Thus, for a fixed point (standard tabloid) $t$ of $D$, define
$$
d^+\! (t) = \sum_{i < j} d_{ij}(t) .
$$
We then have the

\begin{prop}
Suppose $F = \C$, and $D$ is a northwest diagram with $\leq N$ rows and
$\Dh = D$.  Then the betti numbers
$$
\beta_{2i} = \#\{t \mid d^+\! (t) = i\}, \ \ \ \beta_{2i+1} = 0 ,
$$
and the Poincare polynomial
$$
P(x, \FD) \eqdef \sum_i \beta_i x^i = \sum_t x^{2 d^+\! (t)} ,
$$
where $t$ runs over the standard tabloids of $D$.
$\bullet$
\end{prop}

In fact, our proof shows the above propostion for a broader class of
spaces.
Suppose $D = (C_1, C_2, \ldots)$ is an arbitrary
diagram such that the variety
$$
\mbox{\rm Inc}_D \eqdef \{ (V_C)_{C \in D} \in \Gr(D) \mid C \subset C'
\Rightarrow V_C \subset V_{C'} \}
$$
is smooth.  Then the proposition holds with $\FD$ replaced by
$\mbox{\rm Inc}_D$.

\end{document}